\documentclass[12pt]{article}
\usepackage{mathrsfs}
\usepackage{amsmath}
\usepackage{amssymb}
\usepackage{amsthm}
\usepackage{amsfonts}
\usepackage{amscd}
\usepackage[mathscr]{eucal}
\newcommand{\Z} {{\mathbb  Z}}
\newcommand{\Q}{{\mathbb  Q}}
\newcommand{\F}{{\mathbb  F}}

\newcommand{\R} {{\mathbb R}}

\begin{document}
\parindent  25pt
\baselineskip  10mm
\textwidth  15cm    \textheight  23cm \evensidemargin -0.06cm
\oddsidemargin -0.01cm

\title{ { On Some Additive Properties of multiplicative subsemigroups of
Semirings and Arithmetic Applications I}}
\author{\mbox{}
{ Derong Qiu }
\thanks{ \quad E-mail:
derong@mail.cnu.edu.cn; \ derongqiu@gmail.com } \\
(School of Mathematical Sciences,
 Capital Normal University, \\
Beijing 100048, P.R.China )  }

\date{}
\maketitle
\parindent  24pt
\baselineskip  10mm
\parskip  0pt

\par   \vskip 0.4cm

{\bf Abstract} \quad In this paper, we consider a question of sum-keeping about
a multiplicative subsemigroup and its generator subsets in a semiring, and develop
some elementary (collapse) process of the sum-keeping retraction through subsets
until one minimal generators subset. As an application, we study and analyze
several classical problems in additive number theory on the semiring of non-negative
integers by this algebraic and combinatory idea, and provide new proofs in more simple
and direct way for several classical results in number theory. Some further questions
are also presented and discussed.
\par  \vskip  0.2 cm

{ \bf Keywords: } \ Additive property, multiplicative subsemigroup, semiring, prime
number, composite number, Goldbach conjecture, additive number
theory, combinatorics.
\par  \vskip  0.1 cm

{ \bf 2000 Mathematics Subject Classification: } \ 11P32 (primary),
\ 05D05, 11A41, 11T55, 13A99, 16H05 (Secondary).
\par     \vskip  0.4 cm

\hspace{-0.6cm}{\bf 1 \ Introduction}

\par \vskip 0.2 cm

Let $ (A, +, \cdot ) $ be a commutative semiring, $ (M, \cdot ) $
be a multiplicative subsemigroup of $ (A, \cdot ). $ In this paper, we
mainly consider the following question
\par \vskip 0.2 cm

{\bf Question A. } \ Is there a minimal subset $ M_{0} $ of multiplicative generators
of $ M, $ such that $ M + M = M_{0} + M_{0} ? $ that is, $  2 M = 2 M_{0} ? $ (see the
notations below).
\par \vskip 0.2 cm
Here as usual, $ M_{0} $ generates $ M $ means that, for every $ a \in M, $ we have
$ a = b_{1} \cdots b_{n} $ for finitely many elements $ b_{1}, \cdots , b_{n}
\in M_{0}. $ Such $ M_{0} $ is minimal means that any proper subset $ M_{1} $
of $ M_{0} $ can not generates $ M. $
\par \vskip 0.2 cm

For example, in the semiring $ (\Z_{\geq 0}, +, \cdot ) $ of non-negative integers,
if the Goldbach conjecture is true, then the multiplicative subsemigroup $ M =
\Z_{\geq 3 }^{o} $ and its minimal subset of generators $ M_{0} = \mathcal{P}_{\geq 3} $
gives an illustration of Question A on $ \Z_{\geq 0} $ (see Section 3 below for precise
statement and notations). \\
One more general question is
\par \vskip 0.2 cm

{\bf Question B. } \ Let $ (A, +, \cdot ) $ be a commutative semiring, $ \emptyset \neq
A_{0} \subset A, $ find a minimal non-empty subset $ B_{0} \subset A $ such that \\
(1) \ for every $ a \in A_{0}, \ a = b_{1} \cdots b_{n} $ for some $ b_{1}, \cdots , b_{n}
\in B_{0}. $ \\
(2) \ $ A_{0} + A_{0} = B_{0} + B_{0} $ \ (i.e., $ 2 A_{0} = 2 B_{0} $).  \\
Then how to characterize such $ B_{0} $ with $ A_{0} $ and $ A ? $
\par \vskip 0.2 cm

When $ A = \Z_{\geq 0} $ is the semiring of non-negative integers, as a result of
observation and illustration, we have suggested the following
(see the first version [Q, 2009], Conj.3.4)
\par \vskip 0.2 cm

{\bf Conjecture C} (strong form) (See Conjecture 3.4 below and its weak form) \ For every
prime $ p \geq 3, $ one has \ $ 2 \mathcal{S}_{\geq p} = 2 \mathcal{P}_{\geq p}. $ \\
Here $ \mathcal{S}_{\geq p} $ is the multiplicative sub-semigroup
generated by $ \mathcal{P}_{\geq p} $ in $ \Z_{\geq 0} $ (see Def.3.2 below). \\
In other words, Conjecture C says, given any prime $ p \geq 3, $ if an even number $ m $
can be written as $  m = m_{1} + m_{2} $ for two odd numbers $ m_{1}, m_{2} \in \Z_{\geq p}^{o} $
such that no prime factor of $  m_{1} \cdot m_{2} $ is less than $ p, $
then $ m = q_{1} + q_{2} $ for some prime numbers $ q_{1} $ and $ q_{2} $ both of which are equal
or greater than $ p. $ \\
In particular, if we take $ p = 3, $ then since $ \mathcal{S}_{\geq 3} = \Z_{ \geq
3}^{o}, $ this case of Conjecture C above is equivalent to the Goldbach conjecture.
\par \vskip 0.2 cm

The contents of this paper consist of two parts, one is to develop a theory
(in its elementary stage) of the process of sum-keeping retraction via decreasing number
of factors about subsets in a semiring, which is given in section 2 (see e.g., Theorems 2.7,
2.11, 2.15, 2.18 below). For simplicity, we call such process the sum-factor collapse process.
The other one is the arithmetic applications, via the above sum-factor collapse process,
we study some classical problems, mainly on the Goldbach conjecture, in additive
number theory, which is given in section 3 (see e.g., Theorems 3.8,
3.12, 3.15, 3.16 below, see also Remark 3.17 for some explanation of applications).
\par \vskip 0.2 cm

These two parts can read almost independently, and the readers may consult section 2 for the
necessary conceptions and terminologies when they read section 3.
\par \vskip 0.1 cm
The main purpose of this work is try to provide a new way to consider the relation between
addition and multiplication structure in algebra, especially in prime number theory.
In fact, one of my consideration about prime problems is try to find
a sum-keeping process via collapsing all the odd composite numbers from odd integers
(see Example 2.5 and Theorem 3.8 below for the exact meaning). The conception of collapse here
is borrowed from Riemannian geometry I learned around the years 2006 and 2007. This geometric
background motivated the beginning of this work. I wish this elementary idea
would lead to a satisfactory theory to uncover the deep sum-factor
structure in algebra. The other types of sum-factor aspect
and the related applications will be discussed in separated papers.
\par \vskip 0.1 cm

{\bf Notations and terminologies. } \ By a semiring we shall always mean
an associative semiring (which need not be commutative, and need not
contain a multiplicative identity $ 1 $), and we write it as $ (A,
+, \cdot ) $ with addition $ ( + ), $ multiplication $ ( \cdot ), $
and the zero element $ 0 \ ( \neq 1 \ \text{if} \ 1 \in A). $ For $
a \in A, $ we write $ < a > = \{ a, a^{2}, \cdots \} $ be the
multiplicative sub-semigroup generated by $ a. $ Similarly, we write
$ < S > $ for the multiplicative sub-semigroup generated by any
non-empty subset $ S $ in $ A $ (for the definitions, see e.g.,[G],
[Ho]). If $ B_{1}, B_{2} $ and $ C $ are non-empty subsets of $ A, \
b \in A $ and $ n \geq 2 $ is an integer, then
we denote the $ sumsets $ \\
$ B_{1} + B_{2} = \{ b_{1} + b_{2} : \ b_{1} \in B_{1}, \ b_{2} \in
B_{2} \},  \quad b + C = \{ b + c : \ c \in C \}. $ We denote
$ B_{1} \setminus B_{2} = \{ b \in B_{1}: \ b
\notin B_{2} \}. $ Also we denote the $ n-th \ iterated \ sumset $ by \\
$ n C = \{ c_{1} + \cdots + c_{n} : \ c_{1}, \cdots , c_{n} \in C \}
$ and denote the $ n-th \ dilation $ by \\
$ n \ast C = \{ n c : \ c \in C \} $ \
(see [TV], [Na 2] for these standard notations). \\
For the ring $ \Z $ of rational integers and any positive real
number $ x \in \R, $ we denote \\
$ \Z_{\geq x } = \{ k \in \Z : \ k \geq x \}; \quad
 \Z_{\geq x }^{o} = \{ \text{ all odd  integers} \geq x \};  \\
\mathcal{C}^{o}_{\geq x} = \{ \text{all odd  composite numbers} \geq
x \} ;  \quad  \mathcal{P}_{\geq x} = \{ \text{all prime numbers}
\geq x \}. $ \ So \\
$ \Z_{\geq 1 }^{o} = \{ \text{all odd positive integers} \} =
\{ 2 k - 1 : \ k \in \Z \ \text{and} \ k \geq 1 \}; \\
\mathcal{C}^{o}_{> 1} = \{ \text{all odd composite numbers}  \} =
\{\widehat{c}_{1}, \widehat{c}_{2},  \cdots , \widehat{c}_{k},
\cdots  \}, \ ( \widehat{c}_{1} < \widehat{c}_{2} < \cdots ); \\
\mathcal{P} = \mathcal{P}_{\geq 2} = \{ \text{ all prime numbers }
\} = \{ \widehat{p}_{1}, \widehat{p}_{2}, \cdots , \widehat{p}_{k},
\cdots \}, \ ( \widehat{p}_{1} < \widehat{p}_{2} < \cdots ), $ where
$ \widehat{c}_{k} $ is the $ k-$th consecutive odd composite number
and $ \widehat{p}_{k} $ is the $ k-$th consecutive prime number,
e.g., \ $ \widehat{c}_{1} = 9, \widehat{c}_{2} = 15, \widehat{c}_{3}
= 21, \cdots ; \ \widehat{p}_{1} = 2, \widehat{p}_{2} = 3,
\widehat{p}_{3} = 5, \cdots . $ \\
Obviously, $ \Z_{\geq 3 }^{o} = \mathcal{C}^{o}_{> 1} \sqcup
\mathcal{P}_{\geq 3} $ \ (the disjoint union). For any positive
integer $ k, $ we write $ \mathcal{C}^{o}_{> 1}(k) =
\{\widehat{c}_{1}, \cdots , \widehat{c}_{k}\} $ for the set of the
former $ k $ number of odd composite numbers. Obviously, \ $
\mathcal{C}^{o}_{> 1}(1) \subset \mathcal{C}^{o}_{> 1}(2) \subset
\cdots \subset \mathcal{C}^{o}_{> 1}(k) \subset \cdots $ \ and \ $
\mathcal{C}^{o}_{> 1} = \bigcup _{k=1}^{\infty }
\mathcal{C}^{o}_{> 1}(k). $ \\
As usual, $ \pi ( x ) = \Sigma _{ p \in \mathcal{P}, \ p \leq x} 1 $
denotes the number of primes $ \leq x. $ \\
For any set $ X, $ we denote its cardinal by $ \sharp X $ and its
power set by $ 2^{X}, $ which consists of all subsets of $ X. $

\par     \vskip  0.3 cm

\hspace{-0.6cm}{\bf 2 \ The sum-factor collapse process in semirings}

\par \vskip 0.2 cm

To begin with, we introduce several relations $ \gtrdot
\rightleftharpoons _{(m, n)}, \gtrdot \Rightarrow _{(m, n)}, \gtrdot
\rightleftharpoons _{n}, \gtrdot \Rightarrow _{n} $ between
non-empty subsets of a semiring as follows:

{\bf Definition 2.1.} \ Let $ B $ and $ C $ be two non-empty subsets
of a semiring $ A $ and $ m, \ n \in \Z_{ \geq 1}. $ \\
(1) \ We define $ B \gtrdot \rightleftharpoons _{(m, n)} C $ if $ m
B = n C $ and the following condition holds: \\
(FP) \ For any $ b \in B, $ there exist elements $ c_{1}, \cdots ,
c_{k} \in C $ with $ k \geq 1 $ such that $ b = c_{1} \cdots c_{k}.
$ \\
(2) \ We define $ B \gtrdot \Rightarrow _{(m, n)} C $ if both $ B
\supset C $ and $ B \gtrdot \rightleftharpoons _{(m, n)} C $ hold. \\
In particular, if $ m = n > 1, $ then we write simply the
corresponding symbols $ \gtrdot \rightleftharpoons _{(n, n)}, \
\gtrdot \Rightarrow _{(n, n)} $ by $ \gtrdot \rightleftharpoons
_{n}, \ \gtrdot \Rightarrow _{n} $ respectively. For example, $ B
\gtrdot \Rightarrow _{2} C $ means $ B \gtrdot \Rightarrow _{(2, 2)}
C, $ that is, $ B \supset C $ and $ 2 B =2 C $
with the property (FP). \\
Note that a ring is obviously a semiring, and these relations can also be
introduced in rings and other similar algebraic systems.
\par \vskip 0.1 cm
Given non-empty subsets $ B, C, D $ in a semiring $ A, $ by definition,
we have $ B \gtrdot \rightleftharpoons _{n} B. $ Moreover, if $ B
\gtrdot \rightleftharpoons _{n} C $ and $ C \gtrdot
\rightleftharpoons _{n} D, $ then $ B \gtrdot \rightleftharpoons
_{n} D. $ So the relation $ \gtrdot \rightleftharpoons _{n} $ among
non-empty subsets of $ A $ is a pre-order (see [DP]). Particularly,
$ \gtrdot \Rightarrow _{n} $ is a partially ordered relation. A
question here is what conditions will be for $ A $ such that the
relation $ \gtrdot \rightleftharpoons _{n} $ is partially ordered ?
\par \vskip 0.2 cm

{\bf Example 2.2. } \ (1) \ Let $ A = \Z, \ B = \Z_{\geq 1}, \  C=
\{ 2 \} \bigcup \Z_{ \geq 1}^{o}. $ Then obviously $ B \gtrdot
\Rightarrow _{2} C. $ \\
(2) \ Let $ A = \Z / 4 \Z, \ B = \{ \overline{1} \}, \ C = \{
\overline{3} \}. $ Then obviously $ B \gtrdot \rightleftharpoons
_{2} C. $ Moreover, $ 3 B \neq 3 C, $ so $ B \gtrdot
\rightleftharpoons _{3} C $ does not hold. In particular, this
example shows that $  B \gtrdot \rightleftharpoons _{2} C $ does not
imply $ B \supset C $ in general. \\
(3) \ Let $ A $ be an arbitrary ring and $ m, n \in \Z_{\geq 1}, $
we may define a relation $ \geq _{(m, n)} $ between elements in $ A
$ as follows: \ For $ a, b \in A, $ we define $ a \geq _{(m, n)} b $
if $ < a > \gtrdot \rightleftharpoons _{(m, n)} < b >. $  We simply
write $ a \geq _{n} b $ for $ a \geq _{(n, n)} b $ in case $ m = n
> 1. $ Obviously, $ \geq _{n} $ is a pre-order on $ A. $ Let
$ \mathbb{M}_{sg}(A) = \{ < a > : \ a \in A \} $ be the set
consisting of all the multiplicative monogenic sub-semigroups of $ A
$ (see [Ho], p. 8). Then it is easy to see that $ (
\mathbb{M}_{sg}(A), \gtrdot \rightleftharpoons _{n}) $ is a
partially ordered set for every $ n \in \Z_{\geq 1}. $
\par \vskip 0.2 cm

{\bf Lemma 2.3.} \ Let $ B $ and $ C $ be two non-empty subsets of a
semiring $ A $ and $ n \in \Z_{\geq 2}. $ If $ B \gtrdot \Rightarrow
_{n} C, $ then for any subset $ D $ of $ A $ with $ B \supset D
\supset C, $ we have $ B \gtrdot \Rightarrow _{n} D $ and $ D
\gtrdot \Rightarrow _{n} C. $
\par \vskip 0.1 cm

{\bf Proof. } \ It follows easily by the definition.
\par \vskip 0.2 cm

{\bf Definition 2.4.} \ Let $ A $ be a semiring and $ B $ be a non-empty
subset of $ A. $
\par \vskip 0.1 cm
(1) \ For an integer $ n \geq 2, $ we denote \\
$ \mathcal{R}_{d}(2^{B})_{n} = \{ C : \ C \ \text{is a non-empty
subset of} \ B \ \text{and} \ B \gtrdot \Rightarrow _{n} C \}. $ \\
Obviously, $ B \in \mathcal{R}_{d}(2^{B})_{n}. $ For a subset $ C $
of $ B, $ if $ B \setminus C \in \mathcal{R}_{d}(2^{B})_{n}, $ then
we call $ C $ a $ n-$th sum-factor collapsed subset of $ B. $ \\
In particular, if $ \mathcal{R}_{d}(2^{B})_{n} \neq \{ B \} $ \
(respectively, $ \mathcal{R}_{d}(2^{B})_{n} = \{ B \} $), then we
call that $ B $ is $ n$-th sum-factor reducible (respectively,
irreducible). We set \\
$ \mathcal{F}_{r} (A)_{n} = \{ \text{all} \ n-\text{th sum-factor
reducible subsets of} \ A \}, $ \ and \\
$ \mathcal{F}_{ir} (A)_{n} = \{ \text{all} \ n-\text{th sum-factor
irreducible subsets of} \ A \}. $ \\
By definition, if $ \sharp B = 1, $ then $ B \in \mathcal{F}_{ir}
(A)_{n}. $
\par \vskip 0.1 cm
(2) \ Dually, for an integer $ n \geq 2, $ we denote \\
$ \mathcal{E}_{p}(2^{A \supset B})_{n} = \{ D : \ D \ \text{is a
subset of} \ A \ \text{and} \ D \gtrdot \Rightarrow _{n} B \}. $ \\
Obviously, $ B \in \mathcal{E}_{p} (2^{A \supset B})_{n}. $ For a
subset $ C $ of $ A \setminus B, $ if $ B \cup C \in
\mathcal{E}_{p}(2^{A \supset B})_{n}, $ then we call $ C $ a $ n-$th
sum-factor recoverable set for $ B $ in $ A. $ \\
In particular, if $ \mathcal{E}_{p}(2^{A \supset B})_{n} \neq \{ B
\} $ \ (respectively, $ \mathcal{E}_{p}(2^{A \supset B})_{n} = \{ B
\} $), then we call that $ B $ is $ n$-th sum-factor expandable
(respectively, unexpanded). We set \\
$ \mathcal{F}_{e} (A)_{n} = \{ \text{all} \ n-\text{th sum-factor
expandable subsets of} \ A \}, $ \ and \\
$ \mathcal{F}_{ue} (A)_{n} = \{ \text{all} \ n-\text{th sum-factor
unexpanded subsets of} \ A \}. $ \\
By definition, $ A \in \mathcal{F}_{ue} (A)_{n}. $ Obviously $ 2^{A}
\setminus \emptyset = \mathcal{F}_{r} (A)_{n} \sqcup
\mathcal{F}_{ir}(A)_{n} = \mathcal{F}_{e} (A)_{n} \sqcup
\mathcal{F}_{ue}(A)_{n} $ \ (disjoint unions).
\par \vskip 0.2 cm

{\bf Example 2.5.} \ (1) \ For an integer $ n \geq 2, $ it follows by
definition that $ \mathcal{P}_{\geq p} \in \mathcal{F}_{ir} ( \Z
)_{n} $ for each prime number $ p, $ and $ \Z_{ \geq 2 k - 1 }^{o}
\in \mathcal{F}_{ue} (\Z)_{n} $ for all positive integers $ k $
(Generally, for a subset $ B $ in a semiring $ A, $ if $ B $ is a
multiplicative sub-semigroup, then $ B \in \mathcal{F}_{ue} (A)_{n}
$). So $ \Z $ contains infinitely many $ n$-th sum-factor
irreducible subsets $ \mathcal{P}_{\geq p} $ and $ n$-th sum-factor
unexpanded subsets $ \Z_{ \geq 2 k - 1 }^{o}, $ and they form two
descending chains: \ $ \mathcal{P}_{\geq 2} \supset
\mathcal{P}_{\geq 3} \supset \mathcal{P}_{\geq 5} \supset \cdots  $
\ and \ $ \Z_{ \geq 1 }^{o} \supset \Z_{ \geq 3 }^{o} \supset \Z_{
\geq 5 }^{o} \supset \cdots . $
\par \vskip 0.1 cm
(2) \ Let $ B $ be an arithmetic progression of finite length in $
\Z, $ then $ B \in \mathcal{F}_{ue} (\Z)_{n} $ for all $ n \in \Z_{
\geq 2 }. $ \\
We come to prove this conclusion. Write $ k = \sharp B, $ then by
[Na 2], Theorem 1.6 on p.12, one has $ \sharp (n B) = n k - (n - 1).
$ Let $ D \in \mathcal{E}_{p}(2^{\Z \supset B})_{n}, $ then $ B
\subset D \subset \Z $ and $ D \gtrdot \Rightarrow _{n} B, $ so $ n
D = n B. $ Hence $ \sharp (n D) = \sharp (n B) = n k - (n - 1). $ In
particular, $ D $ is a finite subset of $ \Z. $ Write $ l = \sharp
D, $ then $ l \geq k. $ By [Na 2], Theorem 1.3 on p.8, one has $
\sharp (n D) \geq n l - (n - 1), $ which implies $ k \geq l. $ So $
l = k $ and then $ D = B. $ This shows that $ \mathcal{E}_{p}(2^{\Z
\supset B})_{n} = \{ B \}. $ Therefore by definition, $ B \in
\mathcal{F}_{ue} (\Z)_{n}. $ The proof is completed.
\par \vskip 0.1 cm
(3) \ In $ \Z, $ we have $ \Z_{\geq 3 }^{o} \gtrdot \Rightarrow _{2}
\Z_{\geq 3 }^{o} \setminus \{i, \ j \} $ for any $ i, \ j \in
\mathcal{C}^{o}_{> 1}. $ In other words, any pair $ \{ i , \ j \} $
of odd composite numbers is a $ 2-$th sum-factor collapsed subset of
$ \Z_{\geq 3 }^{o}. $ In particular, so is every odd composite
number.  \\
This conclusion can be verified by the unique factorization property of $ \Z.
$ \quad $ \Box $ \\
Note that in Section 3 we will obtain stronger results than the
conclusion (3) of Example 2.5 by a more advanced method (see Theorem
3.8 and Corollary 3.9).
\par \vskip 0.2 cm

{\bf Lemma 2.6.} \ Let $ B $ be a non-empty subset of a ring $ A $
and $ n \in \Z_{\geq 2}. $ Then \\
(1) \ $ B \in \mathcal{F}_{r} (A)_{n} $ if and only if $ B \gtrdot
\Rightarrow _{n} B \setminus \{ b \} $ for some $ b \in B. $ \\
(2) \ In particular, if $ B $ itself is a group under the
multiplication of $ A, $ let $ e $ be the multiplicative unity of $
B $ and $ b_{0} \in B. $ Then we have \\
(2I) \ $ B \gtrdot \Rightarrow _{2} B \setminus \{ b_{0} \} $ \
implies \ $ B \gtrdot \Rightarrow _{2} B \setminus \{ e \} $ \
implies \ $ B \gtrdot \Rightarrow _{2} B \setminus \{ b \} $ for all
$ b \in B. $ \\
(2II) \ $ B \in \mathcal{F}_{r} (A)_{2} $ if and only if $ B \gtrdot
\Rightarrow _{2} B \setminus \{ e \}. $ Particularly, $ B \in
\mathcal{F}_{ir} (A)_{2} $ if $ \sharp B \leq 2. $
\par \vskip 0.1 cm
The conclusion for general $ n $ in case (2) of Lemma 2.6 can be
similarly obtained.
\par \vskip 0.1 cm
{\bf Proof.} \ (1) \ The sufficiency is obvious. For the necessity,
assume $ B \in \mathcal{F}_{r} (A)_{n}, $ then $ B \gtrdot
\Rightarrow _{n} C $ for some non-empty subset $ C $ of $ B $ with $
C \neq B. $ Choose an element $ b \in B \setminus C $ and denote $ D
= B \setminus \{ b \}, $ then $ B \supset D \supset C. $ So by Lemma
2.3, we get $ B \gtrdot \Rightarrow _{n} D. $ \\
(2I) \ Firstly we prove that $ B \gtrdot \Rightarrow _{2} B
\setminus \{ b_{0} \} $ implies $ B \gtrdot \Rightarrow _{2} B
\setminus \{ e \}. $ If $ b = e $ then we are done. So we assume $
b_{0} \neq e. $ Obviously $ e = b_{0} \cdot b_{0}^{-1} $ is a
decomposition of $ e $ in $ B \setminus \{ e \}, $ where $
b_{0}^{-1} $ is the inverse of $ b_{0} $ in $ B. $ For any $ a \in
B, $ one has $ b_{0} \cdot a \in B, $ so by $ B \gtrdot \Rightarrow
_{2} B \setminus \{ b_{0} \} $ we have $ b_{0} + b_{0} \cdot a = c +
d $ for some $ c, d \in B \setminus \{ b_{0} \}. $ Obviously $
b_{0}^{-1} \cdot c, b_{0}^{-1} \cdot d \in B \setminus \{ e \}, $ so
$ e + a = b_{0}^{-1} \cdot (b_{0} + b_{0} \cdot a) = b_{0}^{-1}
\cdot (c + d) = b_{0}^{-1} \cdot c + b_{0}^{-1} \cdot d \in 2 (B
\setminus \{ e \}). $ This implies $ e + B \subset 2 (B \setminus \{
e \}). $ It then follows by the definition that $ B \gtrdot
\Rightarrow _{2} B \setminus \{ e
\}. $ \\
Next we prove that $ B \gtrdot \Rightarrow _{2} B \setminus \{ e \}
$ implies $ B \gtrdot \Rightarrow _{2} B \setminus \{ b \} $ for all
$ b \in B. $ If $ b = e, $ then we are done, so we may assume that $
b \neq e. $ If $ B = \{ e, b \}, $ then $ e + e = b + b $ and $ e +
b = b + b $ because $ B \gtrdot \Rightarrow _{2} B \setminus \{ e \}
= \{ b \}, $ so $ b = e, $ a contradiction! Therefore $ \sharp B >
2, $ and then there exists an element $ c \in B \setminus \{ e, b
\}, $ so $ c^{-1} \in B $ and we have $ e = c \cdot c^{-1}. $
Obviously, $ b \neq b \cdot c^{-1}, $ so $ b = (b \cdot c^{-1})
\cdot c $ is a decomposition of $ b \in B \setminus \{ b \}. $ Now
for any $ a \in B, $ since $ b^{-1} \cdot a \in B, $ by $ B \gtrdot
\Rightarrow _{2} B \setminus \{ e \} $ we have $ e + b^{-1} \cdot a
= c + d $ for some $ c, d \in B \setminus \{ e \}. $ Thus $ bc, bd
\in B \setminus \{ b \} $ and then $ b + a = b \cdot (e + b^{-1}
\cdot a) = b \cdot (c + d) = bc + bd \in 2 (B \setminus \{ b \}). $
This shows that $ b + B \subset 2 (B \setminus \{ b \}), $ and then
it follows easily by definition that $ B \gtrdot \Rightarrow _{2} B
\setminus \{ b \}. $ This proves (2I). \\
(2II) \ follows easily from (2I). The proof of Lemma 2.6 is
completed. \quad $ \Box $
\par \vskip 0.1 cm
As an example, let $ A = \Q [x] $ be the polynomial ring in one
variable $ x, $ and denote $ A^{\times } = \Q ^{\times } = \Q
\setminus \{ 0 \} $ its multiplicative group. Then it is easy to see
that $ \Q ^{\times } \gtrdot \Rightarrow _{2} \Q ^{\times }
\setminus \{ 1 \}, $ and so $ \Q ^{\times } \in \mathcal{F}_{r}
(A)_{2}. $
\par \vskip 0.2 cm

{\bf Theorem 2.7.} \ Let $ B $ be a non-empty subset of a semiring $
A $ and $ n \in \Z_{\geq 2}. $ Then \\
(1) \ the partially ordered set $ ( \mathcal{E}_{p}(2^{A \supset
B})_{n}, \ \subset ) $ contains maximal elements. \\
(2) \ $ \mathcal{E}_{p}(2^{A \supset B})_{n} \cap
\mathcal{F}_{ue}(A)_{n} \neq \emptyset . $ Moreover, an element $ D
\in \mathcal{E}_{p}(2^{A \supset B})_{n} $ is maximal if and only if
$ D \in \mathcal{E}_{p}(2^{A \supset B})_{n} \cap
\mathcal{F}_{ue}(A)_{n}. $
\par \vskip 0.1 cm

{\bf Proof.} \ (1) \ Let $ \{D_{i} \}_{i \in I} $ be a chain in $ (
\mathcal{E}_{p}(2^{A \supset B})_{n}, \ \subset ), $ and let $ D =
\cup _{i \in I} D_{i}. $ Obviously $ D \supset B. $ We want to show
that $ D \gtrdot \Rightarrow _{n} B. $ To see this, firstly, for any
$ d \in D, $ we have $ d \in D_{i} $ for some $ i \in I. $ Since $
D_{i} \gtrdot \Rightarrow _{n} B, $ we have $ d = b_{1} \cdots b_{k}
$ for some $ b_{1}, \cdots , b_{k} \in B $ with $ k \geq 1. $
Secondly, for any $ d_{1}, \cdots , d_{n} \in D, $ there exists an $
j \in I $ such that $ d_{1}, \cdots , d_{n} \in D_{j}. $ Since $
D_{j} \gtrdot \Rightarrow _{n} B, $ we have $ d_{1} + \cdots + d_{n}
\in n D_{j} = n B. $ This shows that $ n D \subset n B $ and hence $
n D = n B $ because $ D \supset B. $ Therefore by definition, we get
$ D \gtrdot \Rightarrow _{n} B, $ i.e., $ D \in \mathcal{E}_{p}(2^{A
\supset B})_{n}. $ So by Zorn's lemma (see, e.g., [AM]), $ (
\mathcal{E}_{p}(2^{A \supset B})_{n}, \ \subset ) $ contains maximal
elements. This proves (1). \\
(2) \ It is easy to see that each maximal element in $
\mathcal{E}_{p}(2^{A \supset B})_{n} $ is $ n$-th sum-factor
unexpanded. So by the result in (1) the inequality as well as the
necessity follow. For the sufficiency, let $ D \in
 \mathcal{E}_{p}(2^{A \supset B})_{n} $ be $ n$-th sum-factor
unexpanded. Then for any $ C \in \mathcal{E}_{p}(2^{A \supset
B})_{n} $ with $ C \supset D, $ since $ C \gtrdot \Rightarrow _{n}
B, $ by Lemma 2.3, we get $ C \gtrdot \Rightarrow _{n} D. $ So $ C =
D $ because $ D $ is $ n$-th sum-factor unexpanded. Therefore $ D $
is maximal in $ \mathcal{E}_{p}(2^{A \supset B})_{n}. $ This proves
(2), and the proof of Theorem 2.7 is completed. \quad $ \Box $ \\
From Theorem 2.7 we know that, for any non-empty subset $ B $ of
a semiring $ A $ and $ n \in \Z_{\geq 2}, $ there exists a $ n-$th
sum-factor unexpanded subset $ D \in \mathcal{F}_{ue}(A)_{n} $ such
that $ D \gtrdot \Rightarrow _{n} B. $
\par \vskip 0.2 cm

{\bf Proposition 2.8.} \ Let $ B $ be a non-empty subset of a semiring $
A. $ If $ B $ is a finite set, then for each $ n \in \Z_{\geq 2}, $
there exists a $ n$-th sum-factor irreducible subset $ C \in
\mathcal{F}_{ir} (A)_{n} $ such that $ B \gtrdot \Rightarrow _{n} C.
$
\par \vskip 0.1 cm

{\bf Proof.} \ Denote by $ t = \sharp B $ the cardinal of $ B, $
then $ t < \infty . $ If $ B \in \mathcal{F}_{ir} (A)_{n}, $ then we
are done. So we assume that $ B \in \mathcal{F}_{r} (A)_{n}. $ Then
by definition, $ B \gtrdot \Rightarrow _{n} B_{1} $ for some
non-empty subset $ B_{1} $ of $ B $ with $ B_{1} \neq B. $ Write $
t_{1} = \sharp B_{1}, $ then $ 1 \leq t_{1} < t. $ If $ B_{1} \in
\mathcal{F}_{ir} (A)_{n}, $ then we are done. Otherwise, we have $
B_{1} \in \mathcal{F}_{r} (A)_{n}, $ so $ B_{1} \gtrdot \Rightarrow
_{n} B_{2} $ for some non-empty subset $ B_{2} $ with $ B_{2} \neq
B_{1}. $ Write $ t_{2} = \sharp B_{2}, $ then $ 1 \leq t_{2} <
t_{1}. $ Following this way, we can obtain a chain  $ B = B_{0}
\supset B_{1} \supset B_{2} \supset \cdots  $ with $ B_{i} \gtrdot
\Rightarrow _{n} B_{i + 1}, \ t_{i} = \sharp B_{i} $ and $ t = t_{0}
> t_{1} > t_{2} > \cdots . $ Obviously this chain must stop after
finite steps because $ t < \infty . $ Assume it stops at $ r $ so we
get a finite chain $ B = B_{0} \supset B_{1} \supset B_{2} \supset
\cdots \supset B_{r} $ with $ B_{r} \in \mathcal{F}_{ir} (A)_{n}. $
Take $ C = B_{r}, $ since the relation $ \gtrdot \Rightarrow _{n} $
is transitive, we get $ B \gtrdot \Rightarrow _{n} C. $ This proves
Proposition 2.8. \quad $ \Box $
\par \vskip 0.2 cm

{\bf Definition 2.9.} \ Let $ B $ be a non-empty subset of a semiring $
A $ and $ a \in A. $ \\
(1) \ If $ a = b_{1} \cdots b_{k} $ for some $ b_{1}, \cdots , b_{k}
\in B $ with $ k \geq 1, $ then we say that $ a $ has a decomposition in
$ B, $ and $ b_{1} \cdots b_{k} $ is an orderly decomposition of $ a
$ in $ B. $ Otherwise, we say that $ a $ has no decompositions in $ B. $
Two orderly decompositions $ a = b_{1} \cdots b_{k} = c_{1} \cdots
c_{h} $ of $ a $ in $ B $ are the same if $ k = h $ and $ b_{i} =
c_{i} $ for each $ i \in \{1, \cdots , k \}. $ For example, if $ a
\in B $ is an idempotent, i.e., $ a = a^{2}, $ then $ a $ has
infinitely many distinct orderly decompositions in $ B $ like the
forms $ a = a^{k} $ as $ k $ varies. \\
If $ a $ has at most finitely many distinct orderly decompositions
in $ B, $ then we call that $ a $ is of finite decomposition in $ B. $ \\
(2) \ For $ b \in B $ and $ r \in \Z_{\geq 1}, $ if there exist $
b_{1}, \cdots , b_{k} \in B $ and $ r_{1}, \cdots , r_{k} \in
\Z_{\geq 1} $ with $ k \geq 1, $ such that $ a = b_{1}^{r_{1}}
\cdots b_{k}^{r_{k}} $ with $ b = b_{i} $ and $ r = r_{i} $ for some
$ i \in \{1, \cdots , k \}, $ then we say that $ b^{r} $ divides $ a $ in
$ B, $ and denote it by $ b^{r} \mid _{B} a. $ Obviously, $ b^{k}
\mid _{B} a $ for any positive integer $ k \leq r. $ If $ b^{r} \mid
_{B} a $ and $ b^{r+1} $ does not divide $ a $ in $ B, $ then we
denote it by $ b^{r} \| _{B} a, $ and call $ r $ the $ b-$order of $
a $ in $ B. $ We write $ r = v_{b, B} (a). $  If $ b $ does not
divide $ a $ in $ B, $ then we write $ v_{b, B} (a) = 0. $ If $
b^{r} \mid _{B} a $ for all positive integers $ r, $ then we write $
v_{b, B} (a) = + \infty . $ \\
(3) \ we define $ F(a \mid B) = \{ b \in B : \ b \mid _{B} a \}, $ \\
$ O(a \mid B) = \{ r \in \Z_{\geq 1} \cup \{ + \infty \} : \ r =
v_{b, B} (a) \ \text{for some} \ b \in B \}, $ and \\
$ T(a \mid B) = \{ (r_{1}, \cdots , r_{k}) : \ a = b_{1}^{r_{1}}
\cdots b_{k}^{r_{k}} \ \text{for some} \ b_{1}, \cdots , b_{k} \in B
\ \text{and} \ r_{1}, \cdots , r_{k} \in \Z_{\geq 1} \}. $ We call
each $ b \in F(a \mid B) $ (resp., each $ (r_{1}, \cdots , r_{k})
\in T(a \mid B) $) a factor (resp., a type of orderly decomposition)
of $ a $ in $ B. $ Obviously, all these sets may be empty. We denote
$ n(a \mid B) = \sharp F(a \mid B), \quad e(a \mid B) = \sup \{ r :
\ r \in O(a \mid B) \}, $ \ and \\
$ d(a \mid B) = \sup \{ r = \sum _{i = 1}^{k} r_{i} : \ (r_{1},
\cdots , r_{k}) \in T(a \mid B) \}, $ and call them the
factor-number, the exponent and the degree of $ a
$ in $ B, $ respectively. \\
If both the factor-number and the degree of $ a $ in $ B $ are
finite, i.e., $ n(a \mid B) < + \infty $ and $ d(a \mid B) < +
\infty , $ then obviously $ a $ is of finite decomposition in $ B. $
If $ A $ is a commutative ring and $ n(a \mid B) < + \infty , $ then
it is easy to see that $ e(a \mid B) < + \infty $ if and only if $
d(a \mid B) < + \infty . $
\par \vskip 0.2 cm

{\bf Definition 2.10.} \ Let $ B $ be a non-empty subset of a semiring
$ A $ and $ a \in A. $ A partition of $ a $ in $ B $ is a finite
sequence of distinct elements $ b_{1}, \cdots , b_{r} \in B $
together with a sequence of positive integers $ k_{1}, \cdots ,
k_{r} $ such that $ a = k_{1} b_{1} + \cdots + k_{r} b_{r}. $ The
number $ k = k_{1} + \cdots + k_{r} $ is called the length of this
partition, and a partition of length $ k $ is called a $
k-$partition. Two partitions $ k_{1} b_{1} + \cdots + k_{r} b_{r} =
k_{1}^{\prime } b_{1}^{\prime } + \cdots + k_{r^{\prime }}^{\prime }
b_{r^{\prime }}^{\prime } $ are the same if there is a bijection $
\varphi : \ \{ 1, \cdots , r \} \rightarrow \{ 1, \cdots , r^{\prime
} \} $ such that $ b_{i} = b_{\varphi (i)}^{\prime } $ and $ k_{i} =
k_{\varphi (i)}^{\prime }. $ If $ a $ has at most finitely many
distinct partitions in $ B, $ then we call that $ a $ is of finite
partition in $ B. $ If $ a $ has at most finitely many distinct $
k-$partitions in $ B, $ then we call that $ a $ is of finite $
k-$partition in $ B. $ For the classical theory of partitions of
positive integers, see [A]. \\
With these definitions, we give the following weaker result which is
dual to the one in Theorem 2.7.
\par \vskip 0.2 cm

{\bf Theorem 2.11.} \ Let $ B $ be a non-empty subset of a semiring
$ A $ and $ n \in \Z_{\geq 2}. $ If every element of $ B $ is of
finite decomposition in $ B, $ and every element of $ n B $ is of
finite $ n-$partition in $ B, $ then the following
statements hold: \\
(1) \ The partially ordered set $ ( \mathcal{R}_{d}(2^{B})_{n}, \
\supset ) $ contains minimal elements. \\
(2) \ $ \mathcal{R}_{d}(2^{B})_{n} \cap \mathcal{F}_{ir}(A)_{n} \neq
\emptyset . $ Moreover, an element $ C \in
\mathcal{R}_{d}(2^{B})_{n} $ is minimal if and only if $ C \in
\mathcal{R}_{d}(2^{B})_{n} \cap \mathcal{F}_{ir}(A)_{n}. $
\par \vskip 0.1 cm

{\bf Proof.} \ (1) \ Let $ \{ C_{i} \}_{i \in I} $ be a chain in $ (
\mathcal{R}_{d}(2^{B})_{n}, \ \supset ), $ and let $ C = \cap _{i
\in I} C_{i}. $ Obviously $ B \supset C. $ We want to show that $ C
\neq \emptyset $ and $ B \gtrdot \Rightarrow _{n} C. $ To see this,
firstly, for any $ b \in B, $ since $ B \gtrdot \Rightarrow _{n}
C_{i} $ for each $ i \in I, $ we have in each $ C_{i} $ at least one
orderly decomposition $ b = c_{i, 1} \cdots c_{i, k_{i}} $ for some
$ c_{i, 1}, \cdots , c_{i, k_{i}} \in C_{i} $ with $ k_{i} \geq 1. $
By assumption, $ b $ is of finite decomposition in $ B. $ Let $ t $
denote the number of all its distinct orderly decompositions in $ B,
$ then $ 1 \leq t < + \infty , $ and such decompositions can
be listed as follows: \\
$ b = b_{m, 1} \cdots b_{m, r_{m}} $ with $ m \in \{ 1, \cdots , t
\}, $ \ where $ b_{m, i} \in B $ and $ r_{m} \in \Z_{\geq 1}. $ \\
If for every $ m \in \{ 1, \cdots , t \}, $ there exists an index $
i_{m} \in I $ such that $ b = b_{m, 1} \cdots b_{m, r_{m}} $ is not
an orderly decomposition of $ b $ in $ C_{i_{m}}, $ i.e., $
C_{i_{m}} $ does not contain the set $ \{ b_{m, 1}, \cdots , b_{m,
r_{m}} \}. $ Then for the minimal element, say $ C_{i_{m_{0}}} $ in
$ \{ C_{i_{1}}, \cdots , C_{i_{t}} \}, \ 1 \leq m_{0} \leq t, $ one
can easily see that $ b $ has no orderly decompositions in $
C_{i_{m_{0}}}, $ a contradiction! Therefore there exists an element
$ m \in \{ 1, \cdots , t \}, $ such that $ b = b_{m, 1} \cdots b_{m,
r_{m}} $ is an orderly decomposition of $ b $ in $ C_{i} $ for each
$ i \in I, $ so $ b_{m, 1}, \cdots , b_{m, r_{m}} \in C_{i} $ for
all $ i \in I. $ Hence $ b_{m, 1}, \cdots , b_{m, r_{m}} \in \cap
_{i \in I} C_{i} = C, $ which gives an orderly decomposition of $ b
$ in $ C. $ In particular, $ C \neq \emptyset . $ This shows that $
B, C $ satisfy the property (FP) of Definition 2.1. \\
Secondly, for any $ b_{1}, \cdots , b_{n} \in B, $ we denote $ b =
b_{1} + \cdots + b_{n} \in n B. $ Since $ B \gtrdot \Rightarrow _{n}
C_{i} $ for each $ i \in I, $ we have in each $ C_{i} $ at least one
$ n-$partition $ b = c_{i, 1} + \cdots + c_{i, n} $ with $ c_{i, 1},
\cdots , c_{i, n} \in C_{i}. $ By assumption, $ b $ is of finite $
n-$partition in $ B. $ Let $ s $ denote the number of all its
partitions of length $ n $ in $ B, $ then $ 1 \leq s < +
\infty , $ and such decompositions can be listed as follows: \\
$ b = b_{j, 1} + \cdots + b_{j, n} $ with $ j \in \{ 1, \cdots , s
\}, $ \ where $ b_{j, 1}, \cdots , b_{j, n} \in B $
(need not be distinct). \\
If for every $ j \in \{ 1, \cdots , s \}, $ there exists an index $
i_{j} \in I $ such that $ b = b_{j, 1} + \cdots + b_{j, n} $ is not
a $ n-$partition of $ b $ in $ C_{i_{j}}, $ i.e., $ C_{i_{j}} $ does
not contain the set $ \{ b_{j, 1}, \cdots , b_{j, n} \}. $ Then for
the minimal element, say $ C_{i_{j_{0}}} $ in $ \{ C_{i_{1}}, \cdots
, C_{i_{s}} \}, \ 1 \leq j_{0} \leq s, $ one can easily see that $ b
$ has no partitions of length $ n $ in $ C_{i_{j_{0}}}, $ a
contradiction! Therefore there exists an element $ j \in \{ 1,
\cdots , s \}, $ such that $ b = b_{j, 1} + \cdots + b_{j, n} $ is a
partition of length $ n $ of $ b $ in $ C_{i} $ for each $ i \in I,
$ so $ b_{j, 1}, \cdots , b_{j, n} \in C_{i} $ for all $ i \in I. $
Hence $ b_{j, 1}, \cdots , b_{j, n} \in \cap _{i \in I} C_{i} = C, $
which gives a partition of length $ n $ of $ b $ in $ C. $ So $ n B
\subset n C. $ Obviously we have $ n C \subset n B $ because $ C
\subset B. $ Therefore $ n B = n C $ and so $ B \gtrdot \Rightarrow
_{n} C, $ i.e., $ C \in \mathcal{R}_{d}(2^{B})_{n}. $ Thus by Zorn's
Lemma, $ ( \mathcal{R}_{d}(2^{B})_{n}, \ \supset ) $ contains
minimal elements. This proves (1). \\
(2) \ It is easy to see that each minimal element in $
\mathcal{R}_{d}(2^{B})_{n} $ is $ n$-th sum-factor irreducible. So
by the result in (1) the inequality as well as the necessity follow.
For the sufficiency, let $ C \in \mathcal{R}_{d}(2^{B})_{n} $ be $
n$-th sum-factor irreducible. Then for any $ D \in
\mathcal{R}_{d}(2^{B})_{n} $ with $ C \supset D, $ since $ B \gtrdot
\Rightarrow _{n} D, $ by Lemma 2.3, we get $ C \gtrdot \Rightarrow
_{n} D. $ So $ C = D $ because $ C $ is $ n$-th sum-factor
irreducible. Therefore $ C $ is minimal in $
\mathcal{R}_{d}(2^{B})_{n}. $ This proves (2), and the
proof of Theorem 2.11 is completed. \quad $ \Box $ \\
It follows that, for any non-empty subset $ B $
of a ring $ A $ and $ n \in \Z_{\geq 2} $ satisfying the condition of
Theorem 2.11, there exists a $ n-$th
sum-factor irreducible subset $ C \in \mathcal{F}_{ir}(A)_{n} $ such
that $ B \gtrdot \Rightarrow _{n} C. $
\par \vskip 0.2 cm

{\bf Example 2.12.} \ Let $ A = \Z $ be the ring of integers and
$ n \in \Z_{\geq 2}. $ \\
(1) \ Take $ B = \mathcal{P}_{\geq 3}, $ by Theorem 2.7, there
exists a $ n-$th sum-factor unexpanded subset $ D \subset \Z $ such
that $ D \gtrdot \Rightarrow _{n} \mathcal{P}_{\geq 3}. $ Obviously
$ D \subset \Z_{ \geq 3 }^{o}. $ \\
(2) \ Take $ B = \Z_{ \geq 3 }^{o}, $ then it is easy to see that $
A $ and $ B $ satisfy the condition of Theorem 2.11, so there
exists a $ n-$th sum-factor irreducible subset $ C \subset \Z $ such
that $ \Z_{ \geq 3 }^{o} \gtrdot \Rightarrow _{n} C. $
Obviously, $ C \supset \mathcal{P}_{\geq 3}. $ \\
Now it is easy to see that the famous Goldbach conjecture is
equivalent to say that $ \Z_{ \geq 3 }^{o} \gtrdot \Rightarrow _{2}
\mathcal{P}_{\geq 3}, $ i.e., $ \mathcal{P}_{\geq 3} \in
\mathcal{R}_{d}(2^{\Z_{ \geq 3 }^{o}})_{2}, $ or equivalently, $
\Z_{ \geq 3 }^{o} \in \mathcal{E}_{p}(2^{\Z \supset
\mathcal{P}_{\geq 3}})_{2}. $ This will be discussed in the next
section. \quad $ \Box $
\par \vskip 0.2 cm

{\bf Proposition 2.13.} \ Let $ B, C $ be two non-empty subsets of a
semiring $ A $ and $ n \in \Z_{\geq 2}. $ If $ B \gtrdot \Rightarrow
_{n} C, $ then for any semiring $ A^{\prime } $ and semiring homomorphism $
\phi : A \rightarrow  A^{\prime } \ (\phi (1) = 1 $ if both $ A $
and $ A^{\prime } $ contain the identities $ 1 $), we have $ \phi
(B) \gtrdot \Rightarrow _{n} \phi (C). $
\par \vskip 0.1 cm

{\bf Proof.} \ Step 1. For any $ b \in B, $ by assumption, $ b =
c_{1} \cdots c_{k} $ for some $ c_{1}, \cdots , c_{k} \in C $ with $
k \geq 1. $ So $ \phi (b) = \phi (c_{1} \cdots c_{k}) = \phi (c_{1})
\cdots \phi (c_{k}) $ with $ \phi (c_{1}), \cdots , \\ \phi (c_{k})
\in \phi (C). $ \ Step 2. For any $ b_{1}, \cdots , b_{n} \in B, $
by assumption, $ b_{1} + \cdots + b_{n} \in n B = n C, $ so $ b_{1}
+ \cdots + b_{n} = c_{1} + \cdots + c_{n} $ for some $ c_{1}, \cdots
, c_{n} \in C. $ Thus $ \phi (b_{1}) + \cdots + \phi (b_{n}) = \phi
(b_{1} + \cdots + b_{n}) = \phi (c_{1} + \cdots + c_{n}) = \phi
(c_{1}) + \cdots + \phi (c_{n}) \in n \phi (C). $ This shows that $
n \phi (B) \subset n \phi (C). $ Moreover, by assumption, $ B
\supset C, $ so $ \phi (B) \supset \phi (C), $ hence $ n \phi (B)
\supset n \phi (C). $ Therefore $ n \phi (B) = n \phi (C), $ and by
definition, we obtain that $ \phi (B) \gtrdot \Rightarrow _{n} \phi
(C). $ This proves Proposition 2.13. \quad $ \Box $
\par \vskip 0.2 cm

{\bf Proposition 2.14.} \ Let $ B $ and $ C $ be two non-empty
subsets of a semiring $ A. $ If $ B \gtrdot \Rightarrow _{2} C, $ then $
B \gtrdot \Rightarrow _{n} C $ for all integers $ n \geq 3. $
\par \vskip 0.1 cm

{\bf Proof.} \ We only need to show that $ n B = n C. $ \\
Case 1. \ $ n = 2 k $ with $ k \in \Z_{\geq 2}. $ For any $ b_{1},
\cdots , b_{2k} \in B, $ since $ 2 B = 2 C, $ we have $ b_{i} + b_{k
+ i} = c_{i} + c_{k + i} $ with $ c_{i}, c_{k + i} \in C $ for each
$ i \in \{1, \cdots , k \}. $ So $ b_{1} + \cdots + b_{2k} = c_{1} +
\cdots + c_{2k} \in 2 k C. $ This shows that $ 2 k B \subset 2 k C $
and so $ 2 k B = 2 k C. $ \\
Case 2. \ $ n = 2 k + 1 $ with $ k \in \Z_{\geq 1}. $ For any $
b_{1}, \cdots , b_{2k}, b_{2k + 1} \in B, $ by case 1 above, we have
$  b_{1} + \cdots + b_{2k} \in 2 k B = 2 k C. $ So $ b_{1} + \cdots
+ b_{2k} = c_{1} + \cdots + c_{2k} $ with $ c_{1}, \cdots , c_{2k}
\in C. $ Also $ c_{2 k} + b_{2 k + 1} \in 2 B = 2 C $ because $ C
\subset B. $ Hence $ c_{2 k} + b_{2 k + 1} = c_{2 k}^{\prime } +
c_{2 k + 1}^{\prime } $ with $ c_{2 k}^{\prime }, c_{2 k +
1}^{\prime } \in C. $ Therefore, $ b_{1} + \cdots + b_{2k} + b_{2k +
1} = c_{1} + \cdots + c_{2k} + b_{2k + 1} = c_{1} + \cdots + c_{2k -
1} + (c_{2 k} + b_{2 k + 1}) = c_{1} + \cdots + c_{2k - 1} + c_{2
k}^{\prime } + c_{2 k + 1}^{\prime } \in (2 k + 1) C. $ This shows
that $ (2 k + 1) B \subset (2 k + 1) C $ and so $ (2 k + 1) B = (2 k
+ 1) C. $ This proves Proposition 2.14. \quad $ \Box $ \\
The converse of the conclusion of Proposition 2.14 is in general not
true. For example, let $ A = \Z / 8 \Z $ be the ring of residue
classes modulo $ 8, $ and $ B = \{\overline{0}, \overline{2},
\overline{4} \}, \ C = \{\overline{0}, \overline{2} \}. $ Then $ B,
C $ satisfy the property (FP) of Definition 2.1. Moreover, it is
easy to see that $ 2 B = 3 B = 3C = \{\overline{0}, \overline{2},
\overline{4}, \overline{6} \}, $ but $ \ 2 C = \{\overline{0},
\overline{2}, \overline{4} \} \neq 2 B. \ $ So we have $ B \gtrdot
\Rightarrow _{3} C, \ B \gtrdot \Rightarrow _{(1, 2)} C, $ and $ B
\gtrdot \Rightarrow _{(2, 3)} C. $ But $ B \gtrdot \Rightarrow _{2}
C $ does not hold. \\
Let $ ( V, + ) $ be an additive group with a partial order $ \leq .
$ Recall that $ V $ is a partial ordered group if it satisfies that
$ v_{1} + w \leq v_{2} + w $ whenever $ v_{1} \leq v_{2} \ (v_{1},
v_{2}, w \in A). $ We write $ v < w $ if $ v \leq w $ and $ v \neq
w. $
\par \vskip 0.2 cm

{\bf Theorem 2.15.} \ Let $ A $ be a ring with the totally
ordered additive group $ (A, + , \leq ) $ and $ b, d \in A $ with $
0 < d, $ set $ B = \{ b + i d : \ i \in \Z_{\geq 0} \}. $ For $ k
\in \Z_{\geq 1}, $ let $ c_{1}, \cdots , c_{k + 1} \in B $ with $ b
< c_{1} < \cdots < c_{k + 1}. $ Denote $ C = \{c_{1}, \cdots , c_{k
+ 1} \} $ and $ D_{i} = C \setminus \{ c_{i} \} $ for each $ i \in
\{ 1, \cdots ,
k + 1 \}. $ If the following conditions hold: \\
(1) \ $ c_{i} $ has a decomposition in $ B \setminus C $ for each $
i \in \{ 1, \cdots , k + 1 \}. $ \\
(2) \ $ b + c_{k + 1} = c_{i_{0}} + c_{j_{0}} $ for some integers $
i_{0}, j_{0} \leq k $ with $ i_{0} \neq j_{0}. $ \\
(3) \ $ B \gtrdot \Rightarrow _{2} B \setminus D_{i} $ for each $ i
\in \{ 1, \cdots , k + 1 \}. $ \\
Then we have $ B \gtrdot \Rightarrow _{2} B \setminus \{c_{1},
\cdots , c_{k + 1} \}, $ that is, $ \{c_{1}, \cdots , c_{k + 1} \} $
is a $ 2-$th sum-factor collapsed subset of $ B. $
\par \vskip 0.1 cm

{\bf Proof.} \ Step 1. By condition (3) for $ i_{0}, $ we get $ b +
c_{k + 1} = x + y $ with $ x, y \in B \setminus D_{i_{0}}. $ If $ x
\in C, $ then obviously $ x = c_{i_{0}} \ ( \neq c_{k + 1}), $ so $
b + c_{k + 1} = c_{i_{0}} + y. $ Thus by condition (2) we get $ y =
c_{j_{0}} \in D_{i_{0}}, $ a contradiction! This shows that $ x
\notin C. $ Similarly, $ y \notin C. $ Therefore $ b + c_{k + 1} = x
+ y \in 2 ( B \setminus C ). $ Step 2. For any $ a \in B, $ by
definition, $ a = b + i d $ for some $ i \in \Z_{\geq 0}. $ We want
to show that $ a + c_{k + 1} \in 2 ( B \setminus C ). $ If $ i = 0,
$ then we are done. So we may assume that $ i > 0. $ Then $ a + c_{k
+ 1} = b + i d + c_{k + 1} = b + (i d + c_{k + 1}) \in 2 ( B
\setminus C) $ because $ b, i d + c_{k + 1} \notin C. $ This shows
that $ c_{k + 1} + B \subset 2 ( B \setminus C ). $ Step 3. For any
$ a_{1}, a_{2} \in B, $ we want to show that $ a_{1} + a_{2} \in 2 (
B \setminus C ). $ If $ a_{1}, a_{2} \notin C, $ then we are done.
Otherwise, we may as well assume that $ a_{1} \in C. $ Then $ a_{1}
= c_{i} $ for some $ i \in \{ 1, \cdots , k + 1 \}. $ If $ i = k +
1, $ then $ a_{1} = c_{k + 1} $ and so by Step 2 we have done.
Otherwise, $ i < k + 1. $ By condition (3), we have $ B \gtrdot
\Rightarrow _{2} B \setminus D_{k + 1}. $ So $ a_{1} + a_{2} = c_{i}
+ a_{2} \in 2 (B \setminus D_{k + 1}), $ thus $ c_{i} + a_{2} =
s_{1} + s_{2} $ with $ s_{1}, s_{2} \in B \setminus D_{k + 1}. $ If
$ s_{1} \neq c_{k + 1} $ and $ s_{2} \neq c_{k + 1}, $ then $ s_{1},
s_{2} \in B \setminus C $ and we are done. Otherwise, we may assume
that $ s_{1} = c_{k + 1}, $ then by Step 2 above we get $ a_{1} +
a_{2} = s_{1} + s_{2} = c_{k + 1} + s_{2} \in 2 (B \setminus C). $
To sum up, we have obtained that $ 2 B \subset 2 (B \setminus C), $
hence $ 2 B = 2 (B \setminus C). $ Also by condition (1), $ B $ and
$ B \setminus C $ satisfy the property (FP) of Definition 2.1.
Therefore $ B \gtrdot \Rightarrow _{2} B \setminus C. $ This proves
Theorem 2.15. \quad $ \Box $
\par \vskip 0.2 cm

In the following, for a ring $ A, $ if $ 1 \in A, $ then we denote
by $ A^{\times } = \{ a \in A : \ \exists b \in A \ \text{s.t.} \ ab
= ba = 1 \} $ the unit group of $ A. $ If $ 1 \notin A, $ then $
A^{\times } = \emptyset . $ For convenience, we call a non-empty
proper subset $ B $ in a commutative ring $ A $ to be
unit-independent if
$$ ( \forall \ b_{1}, b_{2} \in B) \quad b_{1} \cdot A^{\times } \cap
b_{2} \cdot A^{\times } \neq \emptyset \ \Longrightarrow \ b_{1} =
b_{2}. $$ The definition in the non-commutative case is similar. \\
If particularly $ A $ is a UFD (i.e., unique factorization ring, see
[L], p.111), we denote $ E_{ir}(A) = \{ \pi \in A : \ \pi \ \text{is
an irreducible element} \}. $ \\
In general, given two non-empty subsets $ B $ and $ C $ in an
arbitrary semiring $ A, $ if $ B \gtrdot \Rightarrow _{n} C $ \ ($ n \in
\Z_{ \geq 2} $), then obviously $ B \subset <C>, $ where $ <C> $ is
the multiplicative sub-semigroup generated by $ C $ in $ A $ as
above. Those subsets $ B $ satisfying $ <B> \gtrdot \Rightarrow _{n}
B $ will be interesting in themselves.
\par \vskip 0.2 cm

{\bf Definition 2.16.} \ Let $ B $ be a non-empty subset (resp.,
sub-ring, ideal) of a ring $ A $ and $ n \in \Z_{\geq 2}. $ If $ < B
> \gtrdot \Rightarrow _{n} B $ and
$ B \in \mathcal{F}_{ir}(A)_{n}, $ then we call that $ B $ is a $
n-$type optimal subset (resp., sub-ring, ideal) of $ A. $ We denote
\ $ O _{p}(A)_{n} = \{ \text{all} \ n-\text{type optimal subsets of}
\ A \}, \ W_{op} (A) = \{n \in \Z_{\geq 2} : \ O _{p}(A)_{n} \neq
\emptyset \} $ and $ d_{w} (A) = \sharp W_{op} (A). $ We call $
W_{op} (A) $ and $ d_{w} (A) $ the optimal weight set and optimal
degree of $ A $ respectively. If $ W_{op} (A) = \Z_{\geq 2} $
(resp., $ W_{op} (A) = \emptyset $), then we call $ A $ an optimal
(resp., optimal-free) ring, If $ d_{w} (A) = + \infty , $ then we
call $ A $ an almost optimal ring.
\par \vskip 0.2 cm

{\bf Example 2.17.} (1) \ If the Goldbach conjecture would be true,
i.e., $ \Z_{\geq 3 }^{o} \gtrdot \Rightarrow _{2} \mathcal{P}_{\geq
3} $ (see section 3 for the discussion), then by Proposition 2.14,
one has $ \Z_{\geq 3 }^{o} \gtrdot \Rightarrow _{n}
\mathcal{P}_{\geq 3} $ for all $ n \in \Z_{\geq 2}. $ Obviously $
\mathcal{P}_{\geq 3} \in \mathcal{F}_{ir}(\Z)_{n} $  and $ \Z_{\geq
3 }^{o} = < \mathcal{P}_{\geq 3} > , $ so $ \mathcal{P}_{\geq 3} $
is a $ n-$type optimal subset for every $ n \in \Z_{\geq 2}, $ and
by definition, $ \Z $ would be an optimal ring. \\
(2) \ Let $ A = \Z / 6 \Z, \ B = \{\overline{0}, \overline{2} \}. $
Then $ < B > = \{\overline{0}, \overline{2}, \overline{4} \} $ and $
n < B > = n B = < B > $ for all $ n \in \Z_{\geq 2}. $ So $ < B >
\gtrdot \Rightarrow _{n} B $ and $ B \in \mathcal{F}_{ir}(A)_{n}. $
Therefore $ B $ is a $ n-$type optimal subset of $ A $ for every $ n
\in \Z_{\geq 2} $ and so $ A $ is an optimal ring. \\
(3) \ Let $ S $ be a multiplicative sub-semigroup in a ring $ A. $
If $ S \gtrdot \Rightarrow _{n} B $ and $ B \in
\mathcal{F}_{ir}(A)_{n} $ for some subset $ B $ of $ A $ with $ n
\in \Z_{\geq 2}, $ then one can easily see that $ S = < B > $ and $
B $ is a $ n-$type optimal subset of $ A. $
\par \vskip 0.2 cm

{\bf Theorem 2.18.} \ Let $ A = \F_{q}[x] $ be the polynomial
ring in one variable $ x $ over the finite field $ \F_{q} $ of $
q-$elements, and let \\
$ B = \{\text{all irreducible polynomials of degree} \ \geq 1 \
\text{in} \ A \}. $ \\
If $ q $ is odd, then $ B $ is a $ 3-$type optimal subset of $ A. $
\par \vskip 0.1 cm

{\bf Proof.} \ Since $ A $ is a PID (i.e., principal ideal domain),
$ B \in \mathcal{F}_{ir}(A)_{3}. $ Now for any $ a \in \F_{q}, $ we
have $ a = (x + 1) + (q - 2) x + (x + a - 1) \in 3 B, $ so $ \F_{q}
\subset 3 B. $ Moreover, if $ a \neq 0, $ then one always has $ a =
a_{1} + a_{2} + a_{3} $ for some $ a_{1}, a_{2}, a_{3} \in \F_{q}
\setminus \{ 0 \}, $ so $ a x + b = a_{1} x + a_{2} x + (a_{3} x +
b) \in 3 B. $ Next, take $ f_{1}, f_{2}, f_{3} \in < B >, $ write $
f = f_{1} + f_{2} + f_{3}. $ If $ f \in \F_{q} $ or $ f $ is a
polynomial of degree $ 1, $ then by the above discussion, we have $
f \in 3 B. $ If deg $ f \geq 2, $ denote $ f = a g $ with its
leading coefficient $ a \in \F_{q} \setminus \{ 0 \} $ and the monic
polynomial $ g \in A $ with deg $ g = $ deg $ f. $ So $ g $ is a
positive polynomial in the meaning of [EH], hence by the polynomial
$ 3-$primes theorem (see [EH], Theorem A.1 on p.143), we get $ g \in
3 B, $ i.e., $ g = p_{1} + p_{2} + p_{3} $ with positive irreducible
polynomials $ p_{1}, p_{2}, p_{3} \in B, $ so $ f = a g = a p_{1} +
a p_{2} + a p_{3} \in 3 B. $ This implies $ 3 < B > \subset 3 B $
and so $ 3 < B
> = 3 B. $ Therefore $ < B > \gtrdot \Rightarrow _{3} B $ and $ B $
is a $ 3-$type optimal subset of $ A. $ This proves Theorem 2.18.
\quad $ \Box $
\par \vskip 0.1 cm

If $ A = \F_{2}[x] $ is the polynomial ring in one variable $ x $
over the finite field $ \F_{2} $ of $ 2-$elements and  $ B =
\{\text{all irreducible polynomials of degree} \ \geq 1 \ \text{in}
\ A \}, $ then one can easily verify that $ B $ is not a $ 2-$type
optimal subset of $ A. $
\par \vskip 0.2 cm

{\bf Question 2.19.} (1) \ Let $ B $ be a non-empty subset of a ring
$ A. $ Is it true that if $ < B > \gtrdot \Rightarrow _{(m, n)} B $
and $ B \in \mathcal{F}_{ir}(A)_{(m, n)} $ with $ m, n \in \Z_{\geq
2}, $ then $ m = n $ ? Here $ B \in \mathcal{F}_{ir}(A)_{(m, n)} $
means that if $ B \gtrdot \Rightarrow _{(m, n)} C $ for some subset
$ C $ of $ A, $ then $ C = B. $ \\
(2) \ What conditions will be needed for a UFD A such that there
exists a unit-independent subset $ B \subset E_{ir}(A) $ with $
\sharp B = +\infty $ satisfying $ < B > \gtrdot \Rightarrow _{2} B $
(or $ < B > \gtrdot \Rightarrow _{n} B $ for a given $ n \in
\Z_{\geq 2} $) ? For example, if the Goldbach conjecture would be
true, then $ \Z $ would be such a ring.

\par     \vskip  0.3 cm

\hspace{-0.6cm}{\bf 3 \ Arithmetic applications}

\par \vskip 0.2 cm

The famous Goldbach conjecture says that every even integer $ \geq 6
$ is a sum of two odd primes (see [Ch], [HaR], [Nat], [W]), in other
words, it says that $ 2 \Z_{ \geq 3 }^{o} = 2 \mathcal{P}_{\geq 3}.
$ So by the unique factorization of $ \Z, $ this conjecture can be
equivalently stated via the terminology introduced in the above
section 2 as follows:
\par \vskip 0.2 cm

({\bf Goldbach conjecture}) \ $ \Z_{ \geq 3 }^{o} \gtrdot
\Rightarrow _{2} \mathcal{P}_{\geq 3}. $
\par \vskip 0.2 cm

In the following of this section, we mainly consider the ring $ \Z $ of rational integers,
or equivalently, the semiring $ \Z_{\geq 0} $ of non-negative integers. So the notations
used in section above can be significantly simplified, and in most cases, we rewrite in
the standard way. For example, if $ M \gtrdot\Rightarrow_{n} M_{0}, $ where
$ M = < M_{0}> $ is a multiplicative subsemigroup in the semiring
$ (\Z_{\geq 0}, +, \cdot ) $ generated by a subset $ M_{0}, $ then we can rewrite this
simply as $ n M = n M_{0}. $

\par \vskip 0.2 cm

{\bf Proposition 3.1.} \ Let $ k \in \Z_{\geq 0} $ and $ p \in
\mathcal{P}_{\geq 3}. $ If $ \Z_{ \geq 2 k + 1}^{o} \gtrdot
\rightleftharpoons _{n} \mathcal{P}_{\geq p} $ holds for some
integer $ n \geq 2, $ in other words, if $ n \Z_{ \geq 2 k + 1}^{o} =
n \mathcal{P}_{\geq p}, $ and the property (FP) of Def.2.1 above holds.
Then $ k = 1 $ and $ p = 3. $
\par \vskip 0.1 cm

{\bf Proof.} \ It is obvious since $ \Z $ is a UFD, and
$ 3 q \in \Z_{ \geq 2 k + 1}^{o} $ for each odd prime $ q \geq 2k + 1.  $ \quad $ \Box $
\par \vskip 0.2 cm

{\bf Definition 3.2.} \ Let $ p $ be a prime number and $ m, n \in
\Z_{\geq 1}. $ We denote $ \mathcal{S}_{\geq p} = <
\mathcal{P}_{\geq p} > , $ the multiplicative sub-semigroup
generated by $ \mathcal{P}_{\geq p} $ in the semiring $ (\Z_{\geq 0}, +, \cdot ). $
Obviously, $ \mathcal{S}_{\geq p} \in \mathcal{F}_{ue}(\Z_{\geq 0})_{n} $ for each prime
number $ p $ and integer $ n \geq 2 $ (see Example 2.5 (1) above).
Moreover, it is easy to see that $ \mathcal{S}_{\geq 3} = \Z_{ \geq
3}^{o} $ and $ \mathcal{S}_{\geq p} \subset \Z_{ \geq p}^{o} $ but $
\mathcal{S}_{\geq p} \neq \Z_{ \geq p}^{o} $ for all prime numbers $
p \geq 5 $ (for each given $ p \geq 5, $ one always has $ 3^{k} \in
\Z_{ \geq p}^{o} \setminus \mathcal{S}_{\geq p} $ for sufficiently
large integers $ k $). Obviously, one has the strict chain: \ $
\mathcal{S}_{\geq 3} \supsetneqq \mathcal{S}_{\geq 5} \supsetneqq
\cdots . $
\par \vskip 0.2 cm

{\bf Proposition 3.3.} \ (1) \ Let $ p $ and $ q $ be two prime
numbers. If $ \mathcal{S}_{\geq p} \gtrdot \rightleftharpoons _{n}
\mathcal{P}_{\geq q} $ holds for some integer $ n \geq 2, $ then $ p
= q. $ Particularly, if $ n = 2, $ then $ p
= q \geq 3. $ \\
(2) \ If the Goldbach conjecture is true, then $ \mathcal{P}_{\geq
2} $ is a $ n-$type optimal subset for every integer $ n \geq 3. $
\par \vskip 0.1 cm

{\bf Proof.} \ (1) \ It can be directly verified. \\
(2) \ Obviously, $ \mathcal{P}_{\geq 2} \in \mathcal{F}_{ir}(\Z_{\geq 0})_{n}
$ and $ < \mathcal{P}_{\geq 2} > = \Z_{\geq 2}, $ so by the
definition, we only need to show that $ n \Z_{\geq 2} = n
\mathcal{P}_{\geq 2}. $ We use induction on $ n ( \geq 3). $ For $ n
= 3, $ it is easy to see that $ 3 \Z_{\geq 2} = \Z_{\geq 6}. $ For
any $ m \in \Z_{\geq 6}, $ if $ m $ is even, then $ m = 6 = 2 + 2 +
2 \in 3 \mathcal{P}_{\geq 2} $ or $ m \geq 8 $ and $ m - 2 \geq 6 $
is even, so by the Goldbach conjecture, $ m - 2 = p + q $ for some $
p, q \in \mathcal{P}_{\geq 3}. $ Hence $ m = 2 + p + q \in 3
\mathcal{P}_{\geq 2}. $ This shows that $ 3 \Z_{\geq 2} \subset 3
\mathcal{P}_{\geq 2} $ and so $ 3 \Z_{\geq 2} = 3 \mathcal{P}_{\geq
2}. $ Now we assume that our conclusion holds for $ n (\geq 3), $
and we need to show that it also holds for $ n + 1, $ i.e., $ (n +
1) \Z_{\geq 2} = (n + 1) \mathcal{P}_{\geq 2}. $ To see this, let $
b_{1}, \cdots , b_{n + 1} \in \Z_{\geq 2}, $ and denote $ m = b_{1}
+ \cdots + b_{n + 1} = b_{1} + \cdots + b_{n - 1} + (b_{n} + b_{n +
1} - 2) + 2. $ Write $ b_{n}^{\prime } = b_{n} + b_{n + 1} - 2, $
then obviously $ b_{n}^{\prime } \in \Z_{\geq 2}. $ So by the
inductive hypothesis, we get $ b_{1} + \cdots + b_{n - 1} +
b_{n}^{\prime } \in n \Z_{\geq 2} = n \mathcal{P}_{\geq 2}, $ hence
$ b_{1} + \cdots + b_{n - 1} + b_{n}^{\prime } = q_{1} + \cdots +
q_{n} $ for some $ q_{1}, \cdots , q_{n} \in \mathcal{P}_{\geq 2}. $
Therefore $ m = b_{1} + \cdots + b_{n - 1} + b_{n}^{\prime } + 2 =
q_{1} + \cdots + q_{n} + 2 \in (n + 1) \mathcal{P}_{\geq 2}. $ This
shows that $ (n + 1) \Z_{\geq 2} \subset (n + 1) \mathcal{P}_{\geq
2}, $ and so $ (n + 1) \Z_{\geq 2} = (n + 1) \mathcal{P}_{\geq 2}. $
Therefore, by induction, our conclusion holds for all integers $ n
\geq 3. $ This proves (2), and the proof of Proposition 3.3 is
completed. \quad $ \Box $
\par \vskip 0.1 cm

For even integers $ N (> 0), $ let $ G_{2} (N) = \sum _{k_{1} +
k_{2} = N} \Lambda (k_{1}) \Lambda (k_{2}), $ where $ \Lambda (k) $
is the von Mangoldt function, then a stronger form of the Goldbach
problem (which implies its solution for sufficiently large $ N $)
says \\
{\bf Conjecture} (see [IK], p.444). \ For even integers $ N \geq 4,
$ one has \\
$ G_{2} (N) = \mathfrak{S}_{2} (N) N + O (N (\log N)^{-r}), $ \\
where $ \mathfrak{S}_{2} (N) $ is a positive function and $ r $ is a
positive number. see [IK], formula (19.5) on p.444 for the detail.
\\
From this conjecture, one can easily see that the number of
representations of every sufficiently large even integer as a sum of
two odd primes is large, so especially, one can expect that $ 3 +
\mathcal{P}_{\geq k} \subset 2 \mathcal{P}_{\geq 5} $ for a given
positive integer $ k, $ and then predict that $ < \mathcal{P}_{\geq
5} > \gtrdot \Rightarrow _{2} \mathcal{P}_{\geq 5}, $ i.e., $
\mathcal{S}_{\geq 5} \gtrdot \Rightarrow _{2} \mathcal{P}_{\geq 5}.
$ By this observation, we suggest the following
\par \vskip 0.2 cm

{\bf Conjecture 3.4.} (1) (strong form) \ For every prime $ p \geq
3, $ one has \\
$ \mathcal{S}_{\geq p} \gtrdot \Rightarrow _{2}
\mathcal{P}_{\geq p}, $ that is, $ 2 \mathcal{S}_{\geq p} = 2 \mathcal{P}_{\geq p}. $ \\
(2) (weak form) \ For every prime $ p \geq 3, $ there exists a
positive integer $ k $ such that $ \mathcal{S}_{\geq p} \setminus
\mathcal{C}^{o}_{> 1} (k) \gtrdot \Rightarrow _{2} \mathcal{P}_{\geq
p}, $ that is, $ 2 \mathcal{S}_{\geq p} \setminus
\mathcal{C}^{o}_{> 1} (k) = 2 \mathcal{P}_{\geq p}. $
\par \vskip 0.2 cm

{\bf Remark.} In other words, Conjecture 3.4.(1) says, given any prime
$ p \geq 3, $ if an even number $ m $ can be written as $  m = m_{1} + m_{2} $
for two odd numbers $ m_{1}, m_{2} \in \Z_{\geq p}^{o} $ such that no prime factor
of $ m_{1} \cdot m_{2} $ is less than $ p, $ then $ m = q_{1} + q_{2} $ for some prime
numbers $ q_{1} $ and $ q_{2} $ both of which are equal or greater than $ p. $ \\
In particular, if we take $ p = 3, $ then since $ \mathcal{S}_{\geq 3} = \Z_{ \geq
3}^{o}, $ this case of Conj.3.4.(1) above is equivalent to the Goldbach conjecture.
\par \vskip 0.2 cm

Another question here is, can one find all the subsets $ B \subset
\mathcal{P}_{\geq 2} $ such that $ < B > \gtrdot \Rightarrow _{2} B
$ with $ \sharp B = + \infty $ ? One can consider similar questions
about the relations $ \gtrdot \Rightarrow _{n} $ in $ \Z $ for $ n
\in \Z_{\geq 2}. $ Recently, Green and Tao proved the
following famous theorem \\
{\bf Theorem} (Green and Tao [GT]) \ If $ B $ is a subset of prime
numbers with \\
$ \limsup _{N \rightarrow \infty } \sharp (B \cap \{1, \cdots , N
\}) / \pi (N) > 0, $ then for every integer $ k \geq 1, B $ contains
an arithmetic progression of length $ k. $ \\
So naturally, we ask the following
\par \vskip 0.2 cm

{\bf Question 3.5.} \ Can we have $ < B > \gtrdot \Rightarrow _{2} B
$ for every subset $ B $ satisfying the condition of the theorem of
Green-Tao ?
\par \vskip 0.2 cm
Now we come to work out some $ 2-$th sum-factor collapsed subset of
$ \Z_{ \geq 3 }^{o}. $

\par \vskip 0.2 cm

{\bf Proposition 3.6.} \ Let $ B = \{ b_{i} \}_{i = 1}^{\infty } $
be a subsequence of $ \mathcal{C}^{o}_{> 1} = \{ \widehat{c}_{i}
\}_{i = 1}^{\infty }. $ If $ b_{i + 1} - b_{i} > 2 $ for all $ i
\geq 1, $ then $ 2 \Z_{ \geq 3 }^{o} = 2 \Z_{\geq 3 }^{o} \setminus B, $
and $ B $ is a $ 2-$th sum-factor collapsed subset of $ \Z_{ \geq 3 }^{o}. $
\par \vskip 0.1 cm

{\bf Proof.} \ Denote $ H = \Z_{ \geq 3 }^{o} $ and $ D = \Z_{ \geq
3 }^{o} \setminus B. $ Obviously $ \mathcal{P}_{\geq 3} \subset D, $
so by the unique factorization property of $ \Z, $ we only need to
show that $ 2 H \subset 2 D. $ To see this, for any $ c \in
\mathcal{C}^{o}_{> 1}, $ if $ c \notin B, $ then $ 3 + c
 \in 2 D. $ If $ c \in B, $ then $ c = b_{i} $ for some $ i \geq 1.
$ Since $ b_{1} \geq \widehat{c}_{1} \geq 9, $ we have $ 3 + b_{1} =
5 + (b_{1} - 2) \in 2 D $ because $ b_{1} - 2 \in D. $ For $ i > 1,
$ by assumption, $ b_{i} - b_{i - 1} > 2, $ so $ b_{i - 1} < b_{i} -
2 < b_{i}, $ thus $ b_{i} - 2 \notin B, $ hence $ b_{i} - 2 \in D. $
Therefore $ 3 + c = 3 + b_{i} = 5 + (b_{i} - 2) \in 2 D. $ This
shows that $ 3 + \mathcal{C}^{o}_{> 1} \subset 2 D, $ so $ 3 + \Z_{
\geq 3 }^{o} \subset 2 D $ because $ \mathcal{P}_{\geq 3} \subset D
$ and $ \Z_{ \geq 3 }^{o} = \mathcal{C}^{o}_{> 1} \sqcup
\mathcal{P}_{\geq 3} $ (disjoint union). Therefore $ 2 H = 2 \Z_{
\geq 3 }^{o} = 3 + \Z_{ \geq 3 }^{o} \subset 2 D. $ This proves
Proposition 3.6. \quad $ \Box $ \\
From Proposition 3.6, one can obtain many infinite sets which are $
2-$th sum-factor collapsed subsets of $ \Z_{ \geq 3}^{o}. $ For
example, if we let $ \widehat{p}(k)^{o} = \widehat{p}_{2} \cdots
\widehat{p}_{k} $ be the product of all odd primes $ \leq
\widehat{p}_{k} $ for each integer $ k \geq 3, $ and denote $ S =
\{\widehat{p}(k)^{o} : \ k \in \Z_{ \geq 3} \}, $ then $
\widehat{p}(k + 1)^{o} - \widehat{p}(k)^{o} > 2 $ for each integer $
k \geq 3. $ So by Proposition 3.6, one has $ \Z_{ \geq 3 }^{o}
\gtrdot \Rightarrow _{2} \Z_{ \geq 3 }^{o} \setminus S. $ Another
interesting example is the following sequence which Odoni considered
in 1985 (see [R]): \ $ \varpi _{1} = 2, \varpi _{2} = 3, \varpi _{3}
= 7, \cdots , \varpi _{k + 1} = \varpi _{1} \varpi _{2} \cdots
\varpi _{k} + 1. $ Denote $ W = \{\varpi _{k} : \ k \in \Z_{\geq 3}
\} \cap \mathcal{C}^{o}_{> 1}. $ Obviously, $ W $ satisfies the
condition of Proposition 3.6. so $ \Z_{ \geq 3 }^{o} \gtrdot
\Rightarrow _{2} \Z_{ \geq 3 }^{o} \setminus W. $ Therefore, $ S $
and $ W $ are $ 2-$th sum-factor collapsed subsets of $ \Z_{ \geq
3}^{o}. $
\par \vskip 0.2 cm

{\bf Proposition 3.7.} \ For the ring $ \Z, $ the following
statements are equivalent: \\
(1) \ $ \Z_{\geq 3 }^{o} \gtrdot \Rightarrow _{2} \mathcal{P}_{\geq
3}, $ that is, the Goldbach conjecture is true. \\
(2) \ $ \Z_{\geq 3 }^{o} \gtrdot \Rightarrow _{2} \Z_{ \geq 3 }^{o}
\setminus B $ holds for any finite subset $ B $ of
$ \mathcal{C}^{o}_{> 1}, $ that is, any finite subset consisting of odd composite
numbers is a $ 2-$th sum-factor collapsed subset of $ \Z_{ \geq 3 }^{o}. $ \\
(3) \ $ \Z_{\geq 3 }^{o} \gtrdot \Rightarrow _{2} \Z_{ \geq 3 }^{o}
\setminus \mathcal{C}^{o}_{> 1}(k) $ holds for any positive integer
$ k, $ that is, the set of the former $ k \ (\forall k\geq 1) $
number of odd composite numbers is a $ 2-$th sum-factor collapsed subset of
$ \Z_{ \geq 3 }^{o}. $
\par \vskip 0.1 cm

{\bf Proof.} \ (1) $ \Rightarrow $ (2): \ Since \ $ \Z_{\geq 3 }^{o}
\supset \Z_{\geq 3 }^{o} \setminus B \supset \mathcal{P}_{\geq 3}, $
the conclusion follows from Lemma 2.3. \ (2) $ \Rightarrow $ (3) and
(3) $ \Rightarrow $ (1) follow easily by the definitions. \quad $
\Box $
\par \vskip 0.2 cm

Recall that for positive real number $ x, \ \pi ( x ) $ denotes the
number of primes $ \leq x, $ and for every positive integer $ k, \
\widehat{c}_{k} $ is the $ k-$th odd composite number.
\par \vskip 0.2 cm

{\bf Theorem 3.8.} \ Let $ k_{1}, \cdots , k_{r} \in
\mathcal{C}^{o}_{> 1} $ and $ k_{1} < \cdots < k_{r} \ (r \in
\Z_{\geq 3}). $ If $ \pi ( k_{i}) \geq i + 2 $ for all $ i \in \{1,
\cdots , r \}, $ then $ 2\Z_{\geq 3 }^{o} =
2 \Z_{ \geq 3 }^{o} \setminus \{ k_{1}, \cdots , k_{r} \}, $ and
$ \{ k_{1}, \cdots , k_{r} \} $ is a $ 2-$th sum-factor
collapsed subset of $ \Z_{\geq 3 }^{o}. $
\par \vskip 0.1 cm

{\bf Proof.} \ Denote $ B = \Z_{\geq 3 }^{o} $ and $ D = \Z_{ \geq 3
}^{o} \setminus \{ k_{1}, \cdots , k_{r} \}. $ Since $
\mathcal{P}_{\geq 3} \subset D \subset B, $ by the unique
factorization property of $ \Z, $ we only need to show that $ 2 B
\subset 2D. $ To see this, for each $ i \in \{ 1, \cdots , r \}, $
denote $ s_{i} = \pi (k_{i}) $ and let $ \widehat{p}_{1},
\widehat{p}_{2}, \cdots , \widehat{p}_{s_{i}} $ be all the primes $
\leq k_{i} $ with $ \widehat{p}_{1} < \widehat{p}_{2} < \cdots <
\widehat{p}_{s_{i}}. $ As before, note that $ \widehat{p}_{1} = 2,
\widehat{p}_{2} = 3, \widehat{p}_{3} = 5, $ and so on. By
assumption, $ s_{i} \geq i + 2 \geq 3 $
for each $ i. $ Let $ a, b \in B. $ \\
(1 ) \ If $ a, b \notin \{ k_{1}, \cdots , k_{r} \}, $ then $ a, b
\in D, $ so $ a + b \in 2 D. $ \\
(2) \ Otherwise, we may as well assume that $ a \in \{ k_{1}, \cdots
, k_{r} \}, $ and we need to show that $ a + b \in 2 D $ for each $
b \in B. $ \\
(2A) \ Let $ a = k_{r}. $ Firstly, if $ b = 3, $ then $ a + b =
k_{r} + 3 = \widehat{p}_{2} + k_{r} = \widehat{p}_{3} + d_{3} =
\cdots = \widehat{p}_{s_{r}} + d_{s_{r}} $ with $ d_{i} = (k_{r} -
\widehat{p}_{i}) + 3 \in B $ and $ k_{r} > d_{3} > \cdots >
d_{s_{r}} > 3 $ because $ k_{r} - \widehat{p}_{i} > 0 $ for all $ i
\ (3 \leq i \leq s_{r}). $ Since $ \sharp \{d_{3}, \cdots ,
d_{s_{r}} \} = s_{r} - 2 \geq r > r - 1, $ in $ \{d_{3}, \cdots ,
d_{s_{r}} \}, $ there is at least one element, say $ d_{i}, $ such
that $ d_{i} \notin \{ k_{1}, \cdots , k_{r - 1} \}. $ Obviously, $
d_{i} \in D $ and then $ a + b = \widehat{p}_{i} + d_{i} \in 2 D. $
Secondly, if $ b > 3, $ then $ k_{r} + b -3 > k_{r}, $ so $ k_{r} +
b -3 \in D, $ thus $ a + b = 3 + (k_{r} + b -3) \in 2 D. $ This
shows that $ k_{r} + B \subset 2 D. $ \\
(2B) \ Let $ a = k_{i}, \ 1 \leq i < r. $ Suppose that $ k_{j} + B
\subset 2 D $ for each $ j $ such that $ i < j \leq r. $ We come to
verify that $ k_{i} + B \subset 2 D. $ To see this, for any $ b \in
B, $ if $ b = 3, $ then $ a + b = k_{i} + 3 = \widehat{p}_{2} +
k_{i} = \widehat{p}_{3} + d_{3} = \cdots = \widehat{p}_{s_{i}} +
d_{s_{i}} $ with $ d_{l} = (k_{i} - \widehat{p}_{l}) + 3 \in B $ and
$ k_{i} > d_{3} > \cdots > d_{s_{i}} > 3 $ because $ k_{i} -
\widehat{p}_{l} > 0 $ for each $ l \in \{3, \cdots , s_{i} \}. $
Since $ \sharp \{3, \cdots , s_{i} \} = s_{i} - 2 \geq i > i - 1, $
in $ \{d_{3}, \cdots , d_{s_{i}} \}, $ there is at least one
element, say $ d_{l}, $ such that $ d_{l} \notin \{ k_{1}, \cdots ,
k_{i - 1} \}. $ Then obviously $ d_{l} \in D, $ and so $ a + b =
\widehat{p}_{l} + d_{l} \in 2 D. $ If $ b > 3, $ then $ k_{i} + b -
3 > k_{i}. $ If $ k_{i} + b - 3 \in \{k_{i + 1}, \cdots , k_{r} \},
$ i.e., $ k_{i} + b - 3 = k_{j} $ for some $ j \in \{i + 1, \cdots ,
r \}, $ then by our hypothesis, we have $ a + b = 3 + (k_{i} + b -
3) = 3 + k_{j} \in k_{j} + B \subset 2 D. $ If $ k_{i} + b - 3
\notin \{k_{i + 1}, \cdots , k_{r} \}, $ then $ k_{i} + b - 3 \in D,
$ so $ a + b = 3 + (k_{i} + b - 3) \in 2 D. $ This shows that $
k_{i} + B \subset 2 D. $ Therefore, by induction, $ k_{j} + B
\subset 2 D $ for all $ j \in \{1, \cdots , r \}. $ Together with
(1), we obtain $ 2 B \subset 2 D. $ This proves Theorem 3.8. \quad $
\Box $
\par \vskip 0.2 cm

{\bf Corollary 3.9.} \ (1) \ For any $ a_{1}, \cdots , a_{22} \in
\mathcal{C}^{o}_{> 1} $ with $ a_{1} < \cdots < a_{22}, $ we have $
2 \Z_{\geq 3 }^{o} = 2 \Z_{ \geq 3 }^{o} \setminus \{ a_{1}, \cdots , a_{22} \}, $
and $ \{a_{1}, \cdots , a_{22} \} $ is a $ 2-$th sum-factor collapsed
subset of $ \Z_{\geq 3 }^{o}. $ \\
(2) \ Let $ k \in \Z_{\geq 3} $ and $ a_{1}, \cdots , a_{k - 2} \in
\mathcal{C}^{o}_{> 1}. $ If $ a_{i} > \widehat{p}_{k} $ for each $ i
\in \{1, \cdots , k - 2 \}, $ then $ 2 \Z_{\geq 3 }^{o} =
2 \Z_{ \geq 3 }^{o} \setminus \{ a_{1}, \cdots , a_{k- 2} \}, $ and
$ \{ a_{1}, \cdots , a_{k - 2} \} $ is a $ 2-$th sum-factor collapsed subset
of $ \Z_{\geq 3 }^{o}. $
\par \vskip 0.1 cm

{\bf Proof.} \ (1) \ One can directly verify that $ \pi
(\widehat{c}_{i}) \geq i + 2 $ for all integers $ i : \ 1 \leq i
\leq 22, $ but $ \pi (\widehat{c}_{23}) = \pi (93) = 24 < 23 + 2. $
So for any $ a_{1}, \cdots , a_{22} \in \mathcal{C}^{o}_{> 1} $ with
$ a_{1} < \cdots < a_{22}, $ it is easy to see that $ \pi (a_{i})
\geq \pi (\widehat{c}_{i}) \geq i + 2 $ for each $ i \in \{1, \cdots
, 22 \} $ because $ a_{i} \geq \widehat{c}_{i}. $ Therefore by
Theorem 3.8 above, we have $ 2 \Z_{\geq 3 }^{o} = 2 \Z_{ \geq 3 }^{o}
\setminus \{ a_{1}, \cdots , a_{22} \}. $ This proves (1). \\
(2) \ We may as well assume that $ a_{1} < \cdots < a_{k - 2}. $
Then $ \pi (a_{i}) \geq k \geq i + 2 $ for each $ i \in \{1, \cdots
, k - 2 \} $ because $ a_{i} > \widehat{p}_{k} > \cdots >
\widehat{p}_{2} > \widehat{p}_{1}. $ So by Theorem 3.8, we get $
2 \Z_{\geq 3 }^{o} = 2 \Z_{ \geq 3 }^{o} \setminus \{ a_{1}, \cdots , a_{k - 2} \}. $
This proves (2), and the proof of Corollary 3.9 is completed. \quad $ \Box $
\par \vskip 0.2 cm

Note that by Proposition 2.14, it is easy to see that the
conclusions of Propositions 3.6, 3.7, Theorem 3.8 and Corollary 3.9
about the relation $ \gtrdot \Rightarrow _{2} $ also hold for $
\gtrdot \Rightarrow _{n} $ with every integer $ n \geq 2. $
\par \vskip 0.2 cm

{\bf Proposition 3.10.} \ (1) \ For each $ k \in \Z_{\geq 31}, \ \pi
(\widehat{c}_{k}) \leq k - 1 $ and $ \widehat{c}_{k} \leq 4 k - 3. $ \\
(2) \ For each $ k \in \Z_{\geq 3} \setminus \{4, 7 \}, \ 3 +
\widehat{c}_{k} = \widehat{c}_{i} + \widehat{c}_{j} $ with $
\widehat{c}_{i}, \ \widehat{c}_{j} \in \mathcal{C}^{o}_{> 1}(k - 1)
$ and $ i \neq j $ \ (i.e., $ \widehat{c}_{i} \neq \ \widehat{c}_{j}
$).
\par \vskip 0.1 cm

{\bf Proof.} \ Let $ A = \{i \in \Z : \ 1 \leq i \leq \widehat{c}_{k} \}, \ B
= \{2 i : \ i \in \Z \ \text{and} \ 2 \leq 2 i < \widehat{c}_{k} \}.
$ Then $ A \setminus B = \{1, 3, 5, \cdots , \widehat{c}_{k} \} = \{
1 \} \sqcup \mathcal{C}^{o}_{> 1}(k) \sqcup ( \mathcal{P}_{\geq 3}
\cap A) $ (disjoint union). Obviously, $ \sharp ( \mathcal{P}_{\geq
3} \cap A) = \pi (\widehat{c}_{k}) - 1 $ and $ \sharp B =
(\widehat{c}_{k} - 1) / 2. $ Since $ \sharp (A \setminus B) = \sharp
A - \sharp B = \widehat{c}_{k} - (\widehat{c}_{k} - 1) / 2 =
(\widehat{c}_{k} + 1) / 2, $ we get $ (\widehat{c}_{k} + 1) / 2 = 1
+ \sharp \mathcal{C}^{o}_{> 1}(k) + \sharp ( \mathcal{P}_{\geq 3}
\cap A) = 1 + k + \pi (\widehat{c}_{k}) - 1 = k + \pi
(\widehat{c}_{k}). $ So $ \widehat{c}_{k} = 2 k + 2 \pi
(\widehat{c}_{k}) - 1. $ \\
(1) \ By the well known Chebyshev inequality $$ \frac{x}{\log x} <
\pi (x) < 1.25506 \frac{x}{\log x} \ \ (x \geq 17 ) \quad
(\text{see} \ [Nar], p.117), $$ we get $ 2 k + 2 \pi
(\widehat{c}_{k}) - 1 = \widehat{c}_{k} > ( \pi (\widehat{c}_{k})
\cdot \log \widehat{c}_{k}) / 1.25506 \ \ (\widehat{c}_{k} > 17). $
\ So \ $ k > ( \log \widehat{c}_{k} / 2.51012 - 1) \cdot \pi
(\widehat{c}_{k}) \ (\widehat{c}_{k} > 17). $ \ If $ k \geq 46, $
then $ \widehat{c}_{k} \geq \widehat{c}_{46} = 169, $ so $ \log
\widehat{c}_{k} / 2.51012 - 1
> 1, $ thus $ k > \pi (\widehat{c}_{k}), $ i.e.,
$ \pi (\widehat{c}_{k}) \leq k - 1. $ Also $ \widehat{c}_{k} = 2 k +
2 \pi (c_{k}) - 1 \leq 2 k + 2 (k - 1) - 1 = 4 k - 3. $ If $ 31 \leq
k \leq 45, $ then one can directly verify that $ \pi
(\widehat{c}_{k}) \leq k - 1 $ and $ \widehat{c}_{k} \leq 4 k - 3. $
This proves (1). \\
(2) \ Let $ t_{k} = \pi (\widehat{c}_{k}) - 1 $ and $ 3 =
\widehat{p}_{2} < \widehat{p}_{3} < \cdots < \widehat{p}_{t_{k}} $
be the set of all odd primes $ \leq \widehat{c}_{k}. $ If $ k \geq
31, $ then by the proof of (1) above, we know that $ t_{k} \leq k -
2. $ Obviously we have $ 3 + \widehat{c}_{k} = a_{k - 1} +
\widehat{c}_{k - 1} = \cdots = a_{1} + \widehat{c}_{1} $ with $
a_{1}, a_{2}, \cdots , a_{k - 1} \in A \setminus B, $ and $ 3 < a_{k
- 1} < \cdots < a_{2} < a_{1} < \widehat{c}_{k}. $ Since $ (k - 1) -
(t_{k} - 1) \geq 2, $ it is not difficult to see that, in $ \{a_{1},
\cdots , a_{k - 1}\}, $ there are at least two distinct elements,
say $ a_{i}, \ a_{j}, $ not being primes, i.e., $ a_{i}, \ a_{j} \in
\mathcal{C}^{o}_{> 1}(k - 1). $ Now in the equalities $ 3 +
\widehat{c}_{k} = a_{i} + \widehat{c}_{i} = a_{j} + \widehat{c}_{j},
$ if $ a_{i} = \widehat{c}_{i} $ and $ a_{j} = \widehat{c}_{j}, $
then $ 2 a_{i} = 2 a_{j} $ and so $ a_{i} = a_{j}, $ a
contradiction! So $ a_{i} \neq \widehat{c}_{i} $ or $ a_{j} \neq
\widehat{c}_{j}. $ Therefore we obtain the conclusion for all
integers $ k \geq 31. $ If $ 3 \leq k \leq 30 $ and $ k \notin \{4,
7 \}, $ then one can directly verify that $ 3 + \widehat{c}_{k} =
\widehat{c}_{i} + \widehat{c}_{j} \in 2 \mathcal{C}^{o}_{> 1}(k - 1)
$ with $ i \neq j. $ This proves (2), and the proof of Proposition
3.10 is completed. \quad $ \Box $
\par \vskip 0.2 cm

{\bf Corollary 3.11.} \ (1) \ For any $ k \in \Z_{ \geq 3 }
\setminus \{6, 14, 19 \}, $ we have $ 2 k \in 2 \mathcal{C}^{o}_{>
1} \setminus 2 \ast \mathcal{C}^{o}_{> 1} $ or
$ 3 + \mathcal{P}_{\geq 3}. $ \\
(2) \ Let $ k, \ m \in \Z_{\geq 1}. $ If $ \widehat{c}_{k} \geq
13^{(m + 1)}, $ then $ \pi (\widehat{c}_{k}) < \frac{k}{m} $ and $
\widehat{c}_{k} < (2 + \frac{2}{m}) k - 1. $ \\
(3) \ If $ k \geq 31, $ then $ k < ( \frac{1}{2} \log (4 k) - 1 )
\pi (\widehat{c}_{k}) + 1. $
\par \vskip 0.1 cm

{\bf Proof.} \ (1) \ If $ 2 k - 3 \in \mathcal{P}_{\geq 3}, $ then
we are done. Otherwise, $ 2 k - 3 \in \mathcal{C}^{o}_{> 1}, $ then
by Proposition 3.12.(2), we have $ 2 k = 3 + (2 k - 3) \in 2
\mathcal{C}^{o}_{> 1} \setminus 2 \ast \mathcal{C}^{o}_{> 1}. $
This proves (1). \\
(2) \ From the proof of Proposition 3.12.(1), we have $ k > ( \log
\widehat{c}_{k} / 2.51012 - 1) \cdot \pi (\widehat{c}_{k}) \
(\widehat{c}_{k} \geq 17). $ If $ \widehat{c}_{k} \geq 13^{(m + 1)},
$ then $ ( \log \widehat{c}_{k} / 2.51012 - 1) > m, $ so $ k > m
\cdot \pi (\widehat{c}_{k}), $ which implies $ \pi (\widehat{c}_{k})
< \frac{k}{m} $ and $ \widehat{c}_{k} = 2 k + 2 \pi
(\widehat{c}_{k}) - 1 < (2 + \frac{2}{m}) k - 1. $ This
proves (2). \\
(3) \ Since $ \widehat{c}_{31} = 121 > 17, $ for each integer $ k
\geq 31, $ by Proposition 3.12.(1) and its proof, one has $
\widehat{c}_{k} < \pi (\widehat{c}_{k}) \cdot \log
(\widehat{c}_{k}), $ and $ \widehat{c}_{k} \leq 4 k - 3 < 4 k. $ So
$ k = \frac{1}{2} \widehat{c}_{k} - \pi (\widehat{c}_{k}) +
\frac{1}{2} < ( \frac{1}{2} \log (\widehat{c}_{k}) - 1 ) \pi
(\widehat{c}_{k}) + \frac{1}{2} < ( \frac{1}{2} \log (4 k) - 1 ) \pi
(\widehat{c}_{k}) + 1. $ This proves (3), and the proof of Corollary
3.11 is completed. \quad $ \Box $
\par \vskip 0.2 cm

{\bf Theorem 3.12.} \ For every $ k \in \Z_{\geq 1} $ and $ l \in
\Z_{\geq 7} $ satisfying $ 3 + \widehat{c}_{k + l} \geq 2
\widehat{c}_{k}, $ denote $ C (k, l)_{i} = \mathcal{C}^{o}_{> 1} (k
+ l) \setminus \{ \widehat{c}_{k + i} \} \ (i = 0, 1, \cdots , l). $
If $ 2 \Z_{\geq 3 }^{o} = 2 \Z_{\geq 3 }^{o} \setminus C (k, l)_{i} $
for every $ i \in \{0, 1, \cdots , l \}, $
then $ 2 \Z_{\geq 3 }^{o} = 2 \Z_{\geq 3 }^{o} \setminus \mathcal{C}^{o}_{> 1} (k + l). $
\par \vskip 0.1 cm

{\bf Proof.} \ Firstly, we want to verify that $ 3 + \widehat{c}_{k
+ l} \in 2 \mathcal{P}_{\geq 3}. $ In fact, by Proposition 3.10(2),
we have $ 3 + \widehat{c}_{k + l} = \widehat{c}_{i} +
\widehat{c}_{j} $ for some $ \widehat{c}_{i}, \widehat{c}_{j} \in
\mathcal{C}^{o}_{> 1} (k + l - 1) $ with $ i \neq j. $ If $ i < k $
and $ j < k, $ then $ \widehat{c}_{i} < \widehat{c}_{k} $ and $
\widehat{c}_{j} < \widehat{c}_{k}, $ so $ \widehat{c}_{i} +
\widehat{c}_{j} < 2 \widehat{c}_{k} \leq 3 + \widehat{c}_{k + l}, $
a contradiction! so we may as well assume that $ i \geq k. $ Hence $
i = k + i_{0} $ for some $ i_{0} \in \{0, 1, \cdots , l - 1 \}. $ By
our assumption, $ \Z_{\geq 3 }^{o} \gtrdot \Rightarrow _{2} \Z_{\geq
3 }^{o} \setminus C (k, l)_{i_{0}}, $ so $ 3 + \widehat{c}_{k + l}
\in 2 (\Z_{\geq 3 }^{o} \setminus C (k, l)_{i_{0}}), $ i.e., $ 3 +
\widehat{c}_{k + l} = a + b $ for some $ a, b \in \Z_{\geq 3 }^{o}
\setminus C (k, l)_{i_{0}}. $ We assert that $ a, b \in
\mathcal{P}_{\geq 3}. $ If otherwise, we may as well assume that $ a
\notin \mathcal{P}_{\geq 3}, $ then $ a \in \mathcal{C}^{o}_{> 1} (k
+ l), $ and so $ a = \widehat{c}_{k + i_{0}} = \widehat{c}_{i}. $
Hence $ \widehat{c}_{i} + \widehat{c}_{j} = 3 + \widehat{c}_{k + l}
= a + b = \widehat{c}_{i} + b $ and then $ b = \widehat{c}_{j} \in
\mathcal{C}^{o}_{> 1} (k + l) \setminus \{ \widehat{c}_{k + i_{0}}
\} = C (k, l)_{i_{0}} $ because $ i \neq j, $ a contradiction! This
implies $ a, b \in \mathcal{P}_{\geq 3}, $ and so $ 3 +
\widehat{c}_{k + l} = a + b \in 2 \mathcal{P}_{\geq 3}. $ Next,
denote $ B_{l} = \Z_{\geq 3 }^{o} \setminus C (k, l)_{l}, $ we come
to verify that $ B_{l} \gtrdot \Rightarrow _{2} B_{l} \setminus \{
\widehat{c}_{k + l} \}. $ Obviously, $ B_{l} \setminus \{
\widehat{c}_{k + l} \} = \Z_{\geq 3 }^{o} \setminus
\mathcal{C}^{o}_{> 1} (k + l). $ Since $ \mathcal{P}_{\geq 3}
\subset B_{l} \setminus \{ \widehat{c}_{k + l} \}, $ by the unique
factorization property of $ \Z $ and the definition, we only need to
show that $ \widehat{c}_{k + l} + B_{l} \subset 2 (B_{l} \setminus
\{ \widehat{c}_{k + l} \} ). $ To see this, for any $ b \in B_{l}, $
if $ b = 3, $ then by the above discussion, $ \widehat{c}_{k + l} +
b = \widehat{c}_{k + l} + 3 \in 2 \mathcal{P}_{\geq 3} $ and we are
done. If $ b > 3, $ then $ \widehat{c}_{k + l} + b - 3 \in \Z_{\geq
3 }^{o} \setminus \mathcal{C}^{o}_{> 1} (k + l) = B_{l} \setminus \{
\widehat{c}_{k + l} \}, $ so $ \widehat{c}_{k + l} + b = 3 +
(\widehat{c}_{k + l} + b - 3) \in 2 (B_{l} \setminus \{
\widehat{c}_{k + l} \} ). $ Therefore, $ \widehat{c}_{k + l} + B_{l}
\subset 2 (B_{l} \setminus \{ \widehat{c}_{k + l} \} ). $ This
implies $ B_{l} \gtrdot \Rightarrow _{2} B_{l} \setminus \{
\widehat{c}_{k + l} \}. $ Moreover, by our assumption, we have $
\Z_{\geq 3 }^{o} \gtrdot \Rightarrow _{2} \Z_{\geq 3 }^{o} \setminus
C (k, l)_{l} = B_{l}, $ so by the transitivity of $ \gtrdot
\Rightarrow _{2}, $ we get $ \Z_{\geq 3 }^{o} \gtrdot \Rightarrow
_{2} B_{l} \setminus \{ \widehat{c}_{k + l} \}, $ i.e., $ \Z_{\geq 3
}^{o} \gtrdot \Rightarrow _{2} \Z_{\geq 3 }^{o} \setminus
\mathcal{C}^{o}_{> 1} (k + l). $ The proof of Theorem 3.12 is
completed. \quad $ \Box $
\par \vskip 0.2 cm

{\bf Corollary 3.13.} \ For every $ k \in \Z_{\geq 8}, $ denote $
D_{i} = \mathcal{C}^{o}_{> 1} (k) \setminus \{ \widehat{c}_{i} \} \
(i = 1, \cdots , k). $ If $ 2 \Z_{\geq 3 }^{o} =
2 \Z_{\geq 3 }^{o} \setminus D_{i} $ for every $ i \in \{1,
\cdots , k \}, $ then we have $ 2 \Z_{\geq 3 }^{o} =
2 \Z_{\geq 3 }^{o} \setminus \mathcal{C}^{o}_{> 1} (k). $
\par \vskip 0.1 cm

{\bf Proof.} \ It follows easily from the above Theorem 3.12. \quad
$ \Box $
\par \vskip 0.2 cm

{\bf Definition 3.14.} \ We define a sequence of subsets $
\mathcal{C}^{(i)}_{> 1} \ (i = 1, 2, \cdots ) $ of $
\mathcal{C}^{o}_{> 1} $ inductively as follows: \\
$ \mathcal{C}^{(1)}_{> 1} = \{c \in \mathcal{C}^{o}_{> 1}: \ 3 + c
\in 2 \mathcal{P}_{\geq 3} \}, \quad \mathcal{C}^{(2)}_{> 1} = \{c
\in \mathcal{C}^{o}_{> 1} \setminus  \mathcal{C}^{(1)}_{> 1}: \ 3 +
c \in \mathcal{P}_{\geq 3} + \mathcal{C}^{(1)}_{> 1} \}, $ \\
$ \mathcal{C}^{(3)}_{> 1} = \{c \in ( \mathcal{C}^{o}_{> 1}
\setminus \mathcal{C}^{(1)}_{> 1} ) \setminus \mathcal{C}^{(2)}_{>
1}: \ 3 + c \in \mathcal{P}_{\geq 3} + \mathcal{C}^{(2)}_{> 1} \}
\\
= \{c \in \mathcal{C}^{o}_{> 1} \setminus ( \cup _{i = 1}^{2}
\mathcal{C}^{(i)}_{> 1}): \ 3 + c \in \mathcal{P}_{\geq 3} +
\mathcal{C}^{(2)}_{> 1} \}, \ \cdots , $ \\
$ \mathcal{C}^{(k)}_{> 1} = \{c \in \mathcal{C}^{o}_{> 1} \setminus
( \cup _{i = 1}^{k - 1} \mathcal{C}^{(i)}_{> 1}): \ 3 + c \in
\mathcal{P}_{\geq 3} + \mathcal{C}^{(k - 1)}_{> 1} \} \ (k \in
\Z_{\geq 2}). $ \\
Obviously, $ \cup _{i = 1}^{\infty } \mathcal{C}^{(i)}_{> 1} =
\sqcup _{i = 1}^{\infty } \mathcal{C}^{(i)}_{> 1} $ (the disjoint
union). Moreover, it is easy to see that the Goldbach conjecture is
equivalent to say that $ \mathcal{C}^{(1)}_{> 1} =
\mathcal{C}^{o}_{> 1}, $ that is, $ \mathcal{C}^{(k)}_{> 1} =
\emptyset $ for all $ k \in \Z_{\geq 2}. $ We denote
$ \mathcal{C}^{\dag } = \mathcal{C}^{o}_{> 1} \setminus
\mathcal{C}^{(1)}_{> 1}. $ If $ \mathcal{C}^{\dag } \neq \emptyset , $
then we write $ \mathcal{C}^{\dag } = \{c_{1}^{\dag }, c_{2}^{\dag }, \cdots \} $
with $ c_{1}^{\dag } < c_{2}^{\dag } < \cdots , $ and $  c_{1}^{\dag } $
is the smallest odd composite in $ \mathcal{C}^{\dag }. $ Obviously
$ 6 + c_{1}^{\dag } \in 3 \mathcal{P}_{\geq 3}. $ Let $ \mathcal{C}^{\ddag }
= \{c \in \mathcal{C}^{o}_{> 1} : \ 6 + c \notin 3 \mathcal{P}_{\geq 3}\}. $
If $ \mathcal{C}^{\ddag } \neq \emptyset , $ then we write $ \mathcal{C}^{\ddag }
= \{c_{1}^{\ddag }, c_{2}^{\ddag }, \cdots \} $ with $ c_{1}^{\ddag } <
c_{2}^{\ddag } < \cdots , $ and $ c_{1}^{\ddag } $ is the smallest odd
composite in $ \mathcal{C}^{\ddag }. $ Obviously, if $ \mathcal{C}^{\ddag }
\neq \emptyset , $ then $ \mathcal{C}^{\dag } \neq \emptyset , \
c_{1}^{\dag } < c_{1}^{\ddag } $ and $ 9 + c_{1}^{\ddag } \in
4 \mathcal{P}_{\geq 3}. $
\par \vskip 0.2 cm

{\bf Theorem 3.15.} \ For every $ n \in \Z_{\geq 2}, $ we have $
n \Z_{\geq 3 }^{o} = n \Z_{\geq 3 }^{o} \setminus
( \sqcup _{i = 1}^{n - 1} \mathcal{C}^{(i)}_{> 1} ), $ and
$ \sqcup _{i = 1}^{n - 1} \mathcal{C}^{(i)}_{> 1} $ is a $ n-$th
sum-factor collapsed subset of $ \Z_{\geq 3 }^{o}. $
\par \vskip 0.1 cm

{\bf Proof.} \ Since $ \mathcal{P}_{\geq 3} \subset \Z_{\geq 3 }^{o}
\setminus ( \sqcup _{i = 1}^{n - 1} \mathcal{C}^{(i)}_{> 1} ), $ by
the unique factorization property of $ \Z, $ we only need to show
that $ n \Z_{\geq 3 }^{o} \subset n (\Z_{\geq 3 }^{o} \setminus (
\sqcup _{i = 1}^{n - 1} \mathcal{C}^{(i)}_{> 1} )). $ To see this,
firstly we have the following \\
Assertion. $ 3 (k - 1) + \mathcal{C}^{(k - 1)}_{> 1} \subset k
\mathcal{P}_{\geq 3} \quad ( \forall k \in \Z_{\geq 2}). $ \\
We show this by induction on $ k. $ For $ k = 2, $ let $ c_{1} \in
\mathcal{C}^{(1)}_{> 1}, $ by definition, $ 3 + c_{1} \in 2
\mathcal{P}_{\geq 3}. $ Now assume our assertion holds for $ k (
\geq 2), $ we come to consider $ k + 1. $ Let $ c_{k} \in
\mathcal{C}^{(k)}_{> 1}, $ by definition, $ 3 + c_{k} \in
\mathcal{P}_{\geq 3} + \mathcal{C}^{(k - 1)}_{> 1}, $ so $ 3 + c_{k}
= p_{1} + c_{k - 1} $ for some $ p_{1} \in \mathcal{P}_{\geq 3} $
and $ c_{k - 1} \in \mathcal{C}^{(k - 1)}_{> 1}. $ By the inductive
hypothesis, $ 3 (k - 1) + c_{k - 1} \in k \mathcal{P}_{\geq 3}, $
i.e., $ 3 (k - 1) + c_{k - 1} = q_{1} + \cdots + q_{k} $ with $
q_{1}, \cdots , q_{k} \in \mathcal{P}_{\geq 3}, $ so $ 3 k + c_{k} =
3 (k - 1) + (3 + c_{k} ) = 3 (k - 1) + p_{1} + c_{k - 1} = p_{1} + 3
(k - 1) + c_{k - 1} = p_{1} + q_{1} + \cdots + q_{k} \in (k + 1)
\mathcal{P}_{\geq 3}. $ Therefore by induction, the above assertion
holds for all $ k \in \Z_{\geq 2}. $ It then follows that \\
$ 3 (k - 1) + \sqcup _{i = 1}^{k - 1} \mathcal{C}^{(i)}_{> 1}
\subset k \mathcal{P}_{\geq 3}
\quad ( \forall k \in \Z_{\geq 2}). $ \\
In fact, for any $ c \in \sqcup _{i = 1}^{k - 1}
\mathcal{C}^{(i)}_{> 1}, $ we have $ c = c_{i} \in
\mathcal{C}^{(i)}_{> 1} $ for some $ i \in \{1, \cdots , k - 1 \}. $
So by the above assertion, $ 3 i + c_{i} \in (i + 1)
\mathcal{P}_{\geq 3}, $ i.e., $ 3 i + c_{i} = p_{1} + \cdots + p_{i
+ 1} $ for some $ p_{1}, \cdots , p_{i + 1} \in \mathcal{P}_{\geq
3}, $ so $ 3 (k - 1) + c = 3 (k - i - 1) + 3 i + c_{i} = 3 (k - i -
1) + p_{1} + \cdots + p_{i + 1} \in k \mathcal{P}_{\geq 3} $ because
$ 3 \in \mathcal{P}_{\geq 3}. $ Now we come to verify that $ n
\Z_{\geq 3 }^{o} \subset n ( \Z_{\geq 3 }^{o} \setminus ( \sqcup _{i
= 1}^{n - 1} \mathcal{C}^{(i)}_{> 1} ) ). $ To see this, for any $
c_{1}, \cdots , c_{n} \in \Z_{\geq 3 }^{o}, $ denote $ c = c_{1} +
\cdots + c_{n} - 3 (n - 1), $ then obviously $ c \in \Z_{\geq 3
}^{o}. $ If $ c \notin \sqcup _{i = 1}^{n - 1} \mathcal{C}^{(i)}_{>
1}, $ then $ c_{1} + \cdots + c_{n} = 3 (n - 1) + c \in n (
\Z_{\geq 3 }^{o} \setminus ( \sqcup _{i = 1}^{n - 1}
\mathcal{C}^{(i)}_{> 1} ) ). $ If $ c \in \sqcup _{i = 1}^{n - 1}
\mathcal{C}^{(i)}_{> 1} ), $ then by the above discussion, we have $
c_{1} + \cdots + c_{n} = 3 (n - 1) + c \in n \mathcal{P}_{\geq 3}
\subset n ( \Z_{\geq 3 }^{o} \setminus ( \sqcup _{i = 1}^{n - 1}
\mathcal{C}^{(i)}_{> 1} ) ). $ This implies $ n \Z_{\geq 3 }^{o}
\subset n ( \Z_{\geq 3 }^{o} \setminus ( \sqcup _{i = 1}^{n - 1}
\mathcal{C}^{(i)}_{> 1} ) ), $ and the proof of Theorem 3.15 is
completed. \quad $ \Box $
\par \vskip 0.2 cm

{\bf Theorem 3.16.} \ (1) \ $ \mathcal{C}^{o}_{> 1} = \sqcup _{i =
1}^{\infty } \mathcal{C}^{(i)}_{> 1} $ (the disjoint union).
\par \vskip 0.1 cm
(2) \ There exists a positive integer $ k $ such
that $ \mathcal{C}^{(k + i)}_{> 1} = \emptyset $ for all $ i \in
\Z_{\geq 1}, $ in other words, $ \mathcal{C}^{o}_{> 1} = \sqcup _{i
= 1}^{k} \mathcal{C}^{(i)}_{> 1}. $
\par \vskip 0.1 cm
(3) \ Assume $ \mathcal{C}^{\dag } \neq \emptyset . $ Then
we have \\
(a) \ $ 6 + c_{i}^{\dag } \in 3 \mathcal{P}_{\geq 3} \ $ or
$ \ 9 + c_{i}^{\dag } \in 4 \mathcal{P}_{\geq 3} $ \ if \
$ c_{i}^{\dag } > 4 i + 1. $ \\
(b) \ $ 6 + c_{i}^{\dag } \in 3 \mathcal{P}_{\geq 3} $ \ if \
$ \pi (c_{i}^{\dag }) > i + 1. $ \\
(c) \ $ I = \{i \in \Z_{\geq 1} : \ c_{i}^{\dag } \leq 4 i + 1 \} $
is a finite set.
 \par \vskip 0.1 cm

{\bf Proof.} (1) \ If $ \mathcal{C}^{o}_{> 1} \neq \sqcup _{i =
1}^{\infty } \mathcal{C}^{(i)}_{> 1}, $ then there exists a minimal
integer $ c \in \mathcal{C}^{o}_{> 1} $ such that $ c \notin \sqcup
_{i = 1}^{\infty } \mathcal{C}^{(i)}_{> 1}. $ Denote $ r = \pi (c),
$ then $ r \geq 4 $ because $ c \geq \widehat{c}_{1} = 9. $ As
before, let $ 2 = \widehat{p}_{1} < 3 = \widehat{p}_{2} <
\widehat{p}_{3} < \cdots < \widehat{p}_{r} < c $ be all the primes $
< c. $ Then we have $ 3 + c = \widehat{p}_{2} + c_{2} =
\widehat{p}_{3} + c_{3} = \cdots = \widehat{p}_{r} + c_{r} $ with $
c_{2} = c, c_{3}, \cdots , c_{r} \in \Z_{\geq 3 }^{o} $ and $ c =
c_{2} > c_{3} > \cdots > c_{r} > 3. $ If some $ c_{j} \in
\mathcal{P}_{\geq 3}, $ then $ 3 + c = \widehat{p}_{j} + c_{j} \in 2
\mathcal{P}_{\geq 3} $ and so by definition $ c \in
\mathcal{C}^{(1)}_{> 1}, $ a contradiction ! Hence $ c_{j} \in
\mathcal{C}^{o}_{> 1} $ for every $ j \in \{2, 3, \cdots , r \}. $
If among them there is a $ c_{j} \in \sqcup _{i = 1}^{\infty }
\mathcal{C}^{(i)}_{> 1}, $ then $  c_{j} \in \mathcal{C}^{(i)}_{> 1}
$ for some $ i \in \Z_{\geq 1}. $ So $ 3 + c = \widehat{p}_{j} +
c_{j} \in \mathcal{P}_{\geq 3} + \mathcal{C}^{(i)}_{> 1}, $ by
definition, this implies $ c \in \mathcal{C}^{(i + 1)}_{> 1}, $ a
contradiction! Therefore, $ c_{j} \notin \sqcup _{i = 1}^{\infty }
\mathcal{C}^{(i)}_{> 1} $ for any $ j \in \{2, 3, \cdots , r \}. $
But $ r \geq 4 $ and $ c_{j} < c $ for each $ j >2, $ this
contradicts the minimality of $ c. $ So $ \mathcal{C}^{o}_{> 1} =
\sqcup _{i = 1}^{\infty } \mathcal{C}^{(i)}_{> 1}. $ This proves (1). \\
(2) \ Denote $ C = \mathcal{C}^{o}_{> 1} \setminus (
\mathcal{C}^{(1)}_{> 1} \sqcup \mathcal{C}^{(2)}_{> 1}). $ By (1)
above, we only need to show that $ C $ is a finite set. To see this,
by the famous $ 3-$primes theorem of Vinogradov (see e.g., [Na 1],
p. 212), there exists a positive integer $ N_{0} $ such that $ 6 +
\widehat{c}_{i} \in 3 \mathcal{P}_{\geq 2} $ for all $ i \in
\Z_{\geq N_{0}} $ and the number of such representations can be
large enough, for example, greater than $ 6. $ Then we assert that $
\{ \widehat{c}_{i} : \ i \in \Z_{\geq N_{0}} \} \subset
\mathcal{C}^{(1)}_{> 1} \sqcup \mathcal{C}^{(2)}_{> 1}. $ In fact,
for every $ i \geq N_{0}, $ as above, we may take a representation
of $ 6 + \widehat{c}_{i} $ as $ 6 + \widehat{c}_{i} = q_{1} + q_{2}
+ q_{3} $ with $ q_{1}, q_{2}, q_{3} \in \mathcal{P}_{\geq 3} $
because the number of its representations is greater than $ 6. $
Then $ 3 + \widehat{c}_{i} = q_{1} + (q_{2} + q_{3} - 3). $ Denote $
c = q_{2} + q_{3} - 3, $ obviously $ c \in \Z_{\geq 3 }^{o}. $ If $
\widehat{c}_{i} \in \mathcal{C}^{(1)}_{> 1}, $ then we are done.
Otherwise, $ \widehat{c}_{i} \in \mathcal{C}^{o}_{> 1} \setminus
\mathcal{C}^{(1)}_{> 1}, $ then by definition, one can easily see
that $ c \notin \mathcal{P}_{\geq 3}, $ so $ c = \widehat{c}_{j} $
for some $ j \in \Z_{\geq 1}, $ thus $ q_{2} + q_{3} - 3 = c =
\widehat{c}_{j}, $ and then $ 3 + \widehat{c}_{j} = q_{2} + q_{3}
\in 2 \mathcal{P}_{\geq 3}, $ hence by definition, $ \widehat{c}_{j}
\in \mathcal{C}^{(1)}_{> 1} $ and then $ \widehat{c}_{i} \in
\mathcal{C}^{(2)}_{> 1} $ because $ 3 + \widehat{c}_{i} = q_{1} +
\widehat{c}_{j} \in \mathcal{P}_{\geq 3} + \mathcal{C}^{(1)}_{> 1}.
$ This proves our assertion, and so $ C $ is a finite set. This proves (2). \\
(3) \ (a) \ Let $ c_{i}^{\dag } = \widehat{c}_{k} $ be the $ k-$th odd composite
and denote $ t = \pi (c_{i}^{\dag }), $ then
$ 3 = \widehat{p}_{2} < \widehat{p}_{3} < \cdots < \widehat{p}_{t} <
c_{i}^{\dag }, $ and by the proof of proposition 3.10
above we have $ c_{i}^{\dag } = \widehat{c}_{k} = 2 k + 2 t - 1. $ Obviously
$ k > i. $ Denote $ k - i = r, $ and $ \mathcal{C}^{o}_{> 1}(k) \setminus
\{c_{1}^{\dag }, \cdots , c_{i}^{\dag }\} = \{\widehat{c}_{j_{1}}, \cdots ,
\widehat{c}_{j_{r}}\} \subset \mathcal{C}^{(1)}_{> 1} $ with
$ \widehat{c}_{j_{1}} < \cdots < \widehat{c}_{j_{r}}. $ Then
$ 3 + c_{i}^{\dag } = \widehat{p}_{2} + \widehat{c}_{k} = \widehat{p}_{3} + c_{3}
= \cdots = \widehat{p}_{t} + c_{t} $ with $ c_{3}, \cdots , c_{t} \in
\Z_{\geq 3 }^{o} $ and $ \widehat{c}_{k} > c_{3} > \cdots > c_{t} > 3. $
Since $ c_{i}^{\dag } \notin \mathcal{C}^{(1)}_{> 1}, $ we have
$ c_{k} \notin \mathcal{P}_{\geq 3}, $ i.e., $ c_{k} \in \mathcal{C}^{o}_{> 1} $
for each $ k \in \{3, \cdots , t \}. $ \\
Case I. \ If $ \{c_{3}, \cdots , c_{t}\} \cap \{\widehat{c}_{j_{1}}, \cdots ,
\widehat{c}_{j_{r}}\} \neq \emptyset , $ then we have $ c_{h} = \widehat{c}_{j_{l}} $
for some $ h \in \{3, \cdots , t \} $ and $ l \in \{1, \cdots , r \}, $ so
$ 6 + c_{i}^{\dag } = 3 + \widehat{p}_{h} + c_{h} = \widehat{p}_{h} +
(3 + \widehat{c}_{j_{l}}) \in 3 \mathcal{P}_{\geq 3}. $ \\
Case II. \ If $ \{c_{3}, \cdots , c_{t}\} \cap \{\widehat{c}_{j_{1}}, \cdots ,
\widehat{c}_{j_{r}}\} = \emptyset , $ i.e., $ \{c_{3}, \cdots , c_{t}\} \subset
\{ c_{1}^{\dag }, \cdots , c_{i - 1}^{\dag }\}, $ then
$ 3 + c_{i}^{\dag } = \widehat{p}_{2} + \widehat{c}_{k} = \widehat{p}_{3} + c_{3}
= \cdots = \widehat{p}_{t} + c_{t} = d_{j_{1}} + \widehat{c}_{j_{1}} = \cdots =
d_{j_{r}} + \widehat{c}_{j_{r}} $ with $ c_{i}^{\dag } > d_{j_{1}} > \cdots
> d_{j_{r}} > 3 $ and $ d_{j_{1}}, \cdots , d_{j_{r}} \in \mathcal{C}^{o}_{> 1}. $
If $ \{d_{j_{1}}, \cdots , d_{j_{r}} \} \cap \{\widehat{c}_{j_{1}}, \cdots ,
\widehat{c}_{j_{r}} \} = \emptyset , $ then $ \{d_{j_{1}}, \cdots , d_{j_{r}} \}
\subset \{c_{1}^{\dag }, \cdots , c_{i}^{\dag } \} \setminus
\{c_{3}, \cdots , c_{t}, c_{i}^{\dag } \}. $ Denote $ i - (t - 1) - r = s, $
and write $ \{c_{1}^{\dag }, \cdots , c_{i}^{\dag } \} \setminus
(\{d_{j_{1}}, \cdots , d_{j_{r}} \} \sqcup \{c_{3}, \cdots , c_{t},
c_{i}^{\dag } \}) = \{d_{1}^{\prime }, \cdots , d_{s}^{\prime }\} $ with
$ d_{1}^{\prime } < \cdots < d_{s}^{\prime }, $ then
$ 3 + c_{i}^{\dag } = \widehat{p}_{2} + \widehat{c}_{k} = \widehat{p}_{3} + c_{3}
= \cdots = \widehat{p}_{t} + c_{t} = d_{j_{1}} + \widehat{c}_{j_{1}} = \cdots
= d_{j_{r}} + \widehat{c}_{j_{r}} = d_{1}^{\prime } + d_{1}^{\prime \prime } =
\cdots = d_{s}^{\prime } + d_{s}^{\prime \prime } $ with $ d_{1}^{\prime \prime }
> \cdots > d_{s}^{\prime \prime }. $ Obviously $ \{d_{1}^{\prime \prime },
 \cdots , d_{s}^{\prime \prime } \} = \{d_{1}^{\prime }, \cdots , d_{s}^{\prime }\}. $
Since $ \mathcal{C}^{o}_{> 1}(k) = \{c_{1}^{\dag }, \cdots , c_{i}^{\dag } \}
\sqcup \{\widehat{c}_{j_{1}}, \cdots , \widehat{c}_{j_{r}}\} $ and
$ \{c_{1}^{\dag }, \cdots , c_{i}^{\dag } \} =
\{c_{3}, \cdots , c_{t}, c_{i}^{\dag } \} \sqcup
\{d_{j_{1}}, \cdots , d_{j_{r}} \} \sqcup
\{d_{1}^{\prime }, \cdots , d_{s}^{\prime }\}, $ we have
$ k = i + r $ and $ i = t - 1 + r + s. $ Obviously $ s \geq 0, $ so
$ t + r \leq i + 1, $ and then $ c_{i}^{\dag } = 2k + 2t -1 \leq 4 i + 1, $
a contradiction ! Hence $ \{d_{j_{1}}, \cdots , d_{j_{r}} \} \cap
\{\widehat{c}_{j_{1}}, \cdots , \widehat{c}_{j_{r}} \} \neq \emptyset . $
So we have $ d_{j_{h}} = \widehat{c}_{j_{l}} $ for some $ h, l. $
So $ 9 + c_{i}^{\dag } = 6 + d_{j_{h}} + \widehat{c}_{j_{h}} =
(3 + \widehat{c}_{j_{l}}) + (3 + \widehat{c}_{j_{h}}) \in
4 \mathcal{P}_{\geq 3}. $ This proves (a). \\
(b) \ Denote $ t = \pi (c_{i}^{\dag }), $ we have
$ 3 = \widehat{p}_{2} < \widehat{p}_{3} < \cdots < \widehat{p}_{t} <
c_{i}^{\dag }. $ Assume $ t > i + 1. $ By definition, $ 3 + c_{i}^{\dag } =
\widehat{p}_{2} + c_{i}^{\dag } = \widehat{p}_{3} + c_{3} = \cdots =
\widehat{p}_{t} + c_{t} $ with $ c_{3}, \cdots , c_{t} \in \Z_{\geq 3 }^{o} $
and $ c_{i}^{\dag } > c_{3} > \cdots > c_{t} > 3. $ Since $ c_{i}^{\dag }
\notin \mathcal{C}^{(1)}_{> 1}, $ we have $ c_{k} \notin \mathcal{P}_{\geq 3}, $
i.e., $ c_{k} \in \mathcal{C}^{o}_{> 1} $
for each $ k \in \{3, \cdots , t \}. $ Also we have $ \{ c_{3}, \cdots , c_{t} \}
\nsubseteq \{c_{1}^{\dag }, \cdots , c_{i - 1}^{\dag } \} $ because
$ \sharp \{ c_{3}, \cdots , c_{t} \} = t - 2 > i - 1 = \sharp
\{c_{1}^{\dag }, \cdots , c_{i - 1}^{\dag } \}. $ So there exists a $ k \in
\{3, \cdots , t \} $ such that $ c_{k} \in \mathcal{C}^{o}_{> 1} \setminus
\mathcal{C}^{\dag } = \mathcal{C}^{(1)}_{> 1}, $ so $ 3 + c_{k} \in
2 \mathcal{P}_{\geq 3}, $ and then $ 6 + c_{i}^{\dag } = 3 + (3 + c_{i}^{\dag })
= 3 + (\widehat{p}_{k} + c_{k}) = \widehat{p}_{k} + (3 + c_{k}) \in
3 \mathcal{P}_{\geq 3}. $ This proves (b). \\
(c) \ Suppose that $ I $ is an infinite set, then we have an
infinite sequence $ \{ c_{i_{m}}^{\dag }\}_{m=1}^{\infty } $ with
$ c_{i_{1}}^{\dag } < c_{i_{2}}^{\dag } < \cdots < c_{i_{m}}^{\dag }
< \cdots $ and $ c_{i_{m}}^{\dag } \leq 4 i_{m} + 1 $ for every $ i_{m}. $
Take $ \varepsilon = \frac{1}{5}. $ For every positive integer
$ N_{0} > 20, $ there exists an $ i_{m} $ such that $ c_{i_{m}}^{\dag } > N_{0}. $
Denote $ N = 3+ c_{i_{m}}^{\dag } $ and let $ M = \{2k : \ k \in \Z_{\geq 1} \
\text{and} \ 2k \leq N \}, $ then $ \sharp M = \frac{N}{2}, $ and in $ M, $
there are $ i_{m} $ number of even integers which fail to be sums of two odd
primes, i.e., $ 3 + c_{i_{1}}^{\dag }, \cdots , 3 + c_{i_{m}}^{\dag } = N
\notin 2 \mathcal{P}_{\geq 3}. $
Since $ c_{i_{m}}^{\dag } \leq 4 i_{m} + 1, $  we get
$$ \frac{i_{m} + 1}{N} = \frac{i_{m} + 1}{3+ c_{i_{m}}^{\dag }}
\geq \frac{i_{m} + 1}{4 i_{m} + 4} = \frac{1}{4} > \varepsilon +
\frac{1}{N}, \quad \text{i.e.,} \ i_{m} > \varepsilon N, $$
which contradicts to the well known result
of Corput (see [S], p.123$\sim $124). This proves (c), and the
proof of Theorem 3.16 is completed. \quad $ \Box $
\par \vskip 0.2 cm

{\bf Remark 3.17.} \ (1) \ It follows easily from Theorem 3.16(2)
and Theorem 3.15 that
$ n \Z_{\geq 3}^{o} = n \mathcal{P}_{\geq 3} $ for every
sufficiently large integer $ n $ (written as Theorem 3.17(2)
in early versions [Q] of this paper). Yet it is not new since it can be
deduced easily from known classical results of Goldbach problems.
I thank the anonymous expert for pointing out this to me. Nevertheless,
Our proof here is simple and direct.
\par \vskip 0.1 cm
(2) \ Assume $ \mathcal{C}^{\dag } \neq \emptyset . $ If there exists
a positive integer $ i $ such that $ c_{i}^{\dag } \leq 4 i + 1, $ then
by the proof of Proposition 3.10 above, one can easily show that
$ i \geq \frac{c_{1}^{\dag } - 3}{2}. $ Moreover, if $ i $ is the smallest
positive integer such that $ c_{i}^{\dag } \leq 4 i + 1, $ then
$ c_{i}^{\dag } = 4 i + 1 $ and $ c_{i - 1}^{\dag } = 4 i - 1. $
\par \vskip 0.1 cm
(3) \ To show that the set $ I $ of Theorem 3.16 is empty will give a
proof of the statement that every even integer greater than 1 is the
sum of at most four primes. \\
To see this, assume $ I = \emptyset . $ For any even integer $ N $
greater than 1, we come to show that $ N $ can be expressed as the
sum of at most four primes. We may assume that $ N \geq 12. $ Let
$ c = N - 9, $ then $ c \in \Z_{\geq 3 }^{o}. $ If $ c \in
\mathcal{P}_{\geq 3} $ is a prime or $ c \in \mathcal{C}^{(1)}_{> 1}, $
then we are done. So we may assume that $ c \in \mathcal{C}^{\dag }. $ Then
$ c = c_{i}^{\dag } $ for some positive integer $ i. $ By our assumption,
$ I = \emptyset , $ so we have $ c_{i}^{\dag } > 4 i + 1, $ and then by
Theorem 3.16(3)(a), we obtain that $ 6 + c_{i}^{\dag } \in
3 \mathcal{P}_{\geq 3} $ or $ 9 + c_{i}^{\dag } \in
4 \mathcal{P}_{\geq 3}, $ which implies that $ N = 3 + (6 + c_{i}^{\dag })
= 9 + c_{i}^{\dag } \in 4 \mathcal{P}_{\geq 3}. $ The proof is completed.
\par \vskip 0.1 cm
(4) \ Similarly, by Theorem 3.16(3)(b) above, to show that
$ \pi (c_{i}^{\dag }) > i + 1 \ (\forall \
i \in \Z_{\geq 1 }) $ will also give a proof of the statement that every
even integer greater than 1 is the sum of at most four primes.
\par  \vskip 0.2cm

{\bf Remark 3.18.} \ This paper is a revised version of the early one [Q].

\par  \vskip 0.3 cm

\hspace{-0.8cm} {\bf References }
\begin{description}

\item[[A]] G. E. Andrews, The Theory of Partitions, London:
Addison-Wesley Publishing Company, 1976.

\item[[AM]] M.F. Atiyah, I.G. Macdonald, Introduction to Commutative
Algebra, London: Addison-Wesley Publishing Company, 1969.

\item[[Ch]] J.R. Chen, On the representation of a large even integer
as the sum of a prime and the product of at most two primes. { \it
Sci. Sinica,} 16, 157-176 (1973).

\item[[DP]] B.A. Davey, H.A. Priestley, Introduction to Lattices and
Orders, 2nd Edition, Cambridge: Cambridge University Press, 2002.

\item[[EH]] G.W. Effinger, D.R. Hayes, Additive Number Theory
of Polynomials over a Finite Fields, Oxford: Clarendon Press, 1991.

\item[[G]] J.S.Golan, Semirings and Affine Equations over Them:
Theory and Applications. Boston: Kluwer Academic Publishers, 2003.

\item[[GT]] B. Green, T. Tao, The primes contain arbitrarily long
arithmetical progressions, Annals of Math. 167 (2008), 481-547.

\item[[HaR]] H. Halberstam, H.-E.Richert, Sieve Methods, London:
Academic Press, 1974.

\item[[Ho]] J.M. Howie, An Introduction to Semigroup Theory, London:
Academic Press, 1976.

\item[[IK]] H. Iwaniec, E. Kowalski, Analytic Number Theory, AMS,
Providence, Rhode Island, 2004.

\item[[L]] S. Lang, Algebra, 3rd Edition, GTM 211,
New York: Springer-Verlag, 2002.

\item[[Na 1]] M.B. Nathanson, Additive Number Theory-The
Classical Bases, GTM 164, New York: Springer-Verlag, 1996.

\item[[Na 2]] M.B. Nathanson, Additive Number Theory-Inverse
Problems and the Geometry of Sumsets, GTM 165, New York:
Springer-Verlag, 1996.

\item[[Nar]] W. Narkiewicz, The Development of Prime Number Theory,
New York: Springer-Verlag, 2000.

\item[[Q]] D. R. Qiu, Sum-factor decompositions in rings and
arithmetic applications I, arXiv: 0912.4412 v1, 2009;
0912.4412 v2, 2010.

\item[[R]] P. Ribenboim, The Little Book of Bigger Primes,
2nd Edition, New York: Springer-Verlag, 2004.

\item[[S]] W. Sierpinski, Elementary theory of numbers, Editor:
A. Schinzel, North-Holland, North-Holland Math. Library, 1988.

\item[[T]] T. Tao, Every odd number greater than 1 is the sum of
at most five primes, arXiv: 1201.6656v3, 2012.

\item[[TV]] T. Tao, V. Vu, Additive Combinatorics, 2nd Edition,
Cambridge, New York: Cambridge University Press, 2006.

\item[[W]] Y. Wang (edited), Goldbach conjecture, Singapore: World
Scientific Publishing Co Pte Ltd, 1984.

\end{description}

\end{document}